\documentclass[review,onefignum,onetabnum]{siamart171218}

\usepackage[framed,numbered]{matlab-prettifier}
\usepackage{tikz}
\usepackage{hyperref}
\usepackage{cases}
\usepackage[square,sort,comma,numbers]{natbib}
\usepackage{amsmath}   % Mathe 
\usepackage{mathtools} 
\usepackage{subcaption}
\usepackage[font=small,labelfont=small]{caption}
\usepackage{float}
\usepackage{mathrsfs}
\usepackage[normalem]{ulem}
\usepackage[stable]{footmisc}
\usepackage{graphicx}% Include figure files
\usepackage[framed,numbered]{matlab-prettifier}
\usepackage{tikz}
\usepackage{hyperref}
\usepackage{dcolumn}
\usepackage{todonotes}
\usepackage{lipsum}
\usepackage{amsfonts}
\usepackage{graphicx}
\usepackage{tikz}
\usepackage{mathrsfs}
\usepackage{epstopdf}
\usepackage{algorithmicx}
\usepackage{algpseudocode}

\ifpdf
  \DeclareGraphicsExtensions{.eps,.pdf,.png,.jpg}
\else
  \DeclareGraphicsExtensions{.eps}
\fi

\newsiamremark{remark}{Remark}
\newsiamremark{hypothesis}{Hypothesis}
\crefname{hypothesis}{Hypothesis}{Hypotheses}
\newsiamthm{claim}{Claim}

\headers{ Optimal control of thin films over flexible walls }{Alrashidy et al.}

\title{ Optimal Control of THIN-FILM FLOW on a Flexible Topography \thanks{Submitted to the editors DATE.
\funding{this work was funded by Islamic University of Madinah and Ministry of Education of Saudi Arabia. The research by KvdZ was supported by the Engineering and Physical Sciences Research Council (EPSRC), UK under Grant and EP/W010011/1. }}}

\author{ 
S. Alrashidy\thanks{School of Mathematical Sciences, University of Nottingham, University Park, Nottingham NG7 2RD, UK (\email{sattam.alrashidy@nottingham.ac.uk}). Department of Mathematics, Islamic University of Madinah, Madinah, Saudi Arabia (\email{sattamjr@iu.edu.sa}).}
\and
A. Kalogirou\thanks{School of Mathematical Sciences, University of Nottingham, University Park, Nottingham NG7 2RD, UK (\email{anna.kalogirou@nottingham.ac.uk}).}
\and
D. Kalise\thanks{Department of Mathematics, Imperial College London, London SW7 2AZ, UK (\email{d.kalise-balza@imperial.ac.uk}).}
\and
K.~G. van der Zee\thanks{School of Mathematical Sciences, University of Nottingham, University Park, Nottingham NG7 2RD, UK (\email{kg.vanderzee@nottingham.ac.uk}).}\thanks{Department of Mathematics, Imperial College London, London SW7 2AZ, UK}
}

\usepackage{amsopn}

\makeatletter
\newcommand*{\addFileDependency}[1]{% argument=file name and extension
  \typeout{(#1)}% latexmk will find this if $recorder=0 (however, in that case, it will ignore #1 if it is a .aux or .pdf file etc and it exists! if it doesn't exist, it will appear in the list of dependents regardless)
  \@addtofilelist{#1}% if you want it to appear in \listfiles, not really necessary and latexmk doesn't use this
  \IfFileExists{#1}{}{\typeout{No file #1.}}% latexmk will find this message if #1 doesn't exist (yet)
}
\makeatother

\ifpdf
\hypersetup{
  pdftitle={An Example Article},
  pdfauthor={D. Doe, P. T. Frank, and J. E. Smith}
}
\fi

\usepackage{xcolor}
\usepackage{colortbl} % for \rowcolor
\usepackage[normalem]{ulem} % for \sout
\definecolor{r1}{RGB}{0,0,160}   % Reviewer 1 (blue)
\definecolor{r2}{RGB}{0,110,0}   % Reviewer 2 (green)
\definecolor{gone}{gray}{0.45}   % deleted text (grey)
\usepackage[normalem]{ulem} % for \sout
\usepackage{rotating} % for sidewaystable environment
\usepackage[normalem]{ulem}  % for \sout
\usepackage{xcolor}
\definecolor{r1}{RGB}{0,0,180}     % blue for Reviewer 1
\definecolor{r2}{RGB}{0,130,0}     % green for Reviewer 2
\definecolor{gone}{RGB}{120,120,120} % grey for deleted text

\definecolor{rb}{RGB}{130, 0, 130} % purple-like

\usepackage{url}
\usepackage{hyperref}
\usepackage{natbib}
\nolinenumbers
\begin{document}

\maketitle

%% REQUIRED
\begin{abstract}
   This work presents a mathematical model for the optimal control of thin-film flows over a flexible substrate influenced by an external force. The objective is to find the optimal distributed force acting on the topography that minimises the differences between actual and desired thin-film profiles. 
A nonlinear lubrication equation governing the fluid dynamics and appropriate functional settings for this model are presented. 
It is also shown that this system satisfies a global energy-dissipation law for a suitable energy functional. 
Optimality conditions are derived for the solution of the minimisation problem of a specified cost function across a time horizon. 
These conditions are formulated at a continuous level as system of coupled, forward-backward PDEs, which are subsequently discretised for numerical investigation.
To ensure computational efficiency and stability, first-order Implicit-Explicit (IMEX) time-stepping schemes are employed to handle the nonlinearities in the model, and a reduced gradient descent algorithm is applied to obtain a numerical approximation of the optimal control signal. 
Numerical results illustrate that controlling the thin film, even during rupture, achieves a precise film profile. 
This control strategy accelerates convergence towards a steady state, reduces instabilities, stabilises dewetting processes, and meets the desired profile specifications.

\end{abstract}

% REQUIRED

% REQUIRED

\section{Introduction}

Controlling thin-film flow has gained significant interest across a range of applications~\cite{wray2022electrostatic,kistler1997coating,fowler2016controlling,miyara1999numerical,serifi2004transient}. For example, coating processes favour stable, flat-flow films~\cite{tomlin2019optimal}, whereas applications in heat and mass transfer often require interfacial deformations~\cite{frisk1972enhancement}. Motivation stems from the increasing number of thin-film devices manufactured across diverse industries such as microelectronics, displays, optical storage, and microfluidic devices~\cite{kistler1997coating,sellier2010beating,fowler2016controlling,miyara1999numerical,serifi2004transient}. Of particular interest are inverted film configurations where the liquid resides on the underside of a substrate~\cite{blyth2023transition}. Such configurations arise in liquid film coating in cooling towers~\cite{rohlfs2017hydrodynamic}, environmental applications such as glacier hydrology and cave pattern morphogenesis~\cite{camporeale2017asymptotic}, and biomedical applications such as microbicide epithelial coatings and blood pattern analysis~\cite{kabaliuk2013blood,hu2016effects}. In these inverted geometries, gravity destabilises the film and control of the thin layer deposition process is required to achieve desired free surface profiles~\cite{decre2003gravity}.
%These applications all require a thorough understanding and control of the thin layer deposition process and the free surface profile~\cite{sellier2010beating,decre2003gravity}.
However, efficient control of interfacial shape or surface area has been largely unexplored until recently, with most efforts focusing on varying external mechanisms like electric fields~\cite{anderson2017electric,tseluiko2006wave}, including a rigid topography~\cite{gaskell_jimack_sellier_thompson_wilson_2004,gaskell2004gravity,sellier2016inverse,tseluiko2013stability,sellier2008substrate} or surfactants~\cite{blyth2004effect}, rather than using mathematically optimal or robust approaches~\cite{wray2022electrostatic}. 

%MOVE UP AS SECOND PARAGRAPH, REPRHASE STARTING SENTENCE:
The mathematical model of thin films is based on the lubrication approximation simplified from the Navier-Stokes equations~\cite{reynolds1886iv}, resulting in a time-dependent Partial Differential Equation (PDE) for the height of the film. This model has been extensively studied for its ability to facilitate accurate analytical and numerical calculations %, often employing the long-wave assumption and formal asymptotic~
\cite{alexander2020stability,lee2011dynamics,matar2007falling,matar2004rupture}. This indicates that the characteristic wavelengths of interfacial deformation are significantly larger than the film thickness ratio~\cite{tomlin2019optimal,oron1997long}. Additionally, when combined with flexible substrates, further dynamics arise due to elastohydrodynamics~\cite{matar2007falling,matar2004rupture}. Craster and Matar's review provides a comprehensive overview of these models~\cite{craster2009dynamics}. Traditional models typically assume rigid substrates, but real-world applications increasingly involve flexible substrates~\cite{matar2004rupture}.

%MAKE MORE SCIENTIFIC/NEUTRAL: 

Our work develops a PDE-constrained optimisation approach for the control of thin-film flows over a flexible substrate. Our thin film model allows for the rupture of films as well as the coalescence of bubbles. We focus on finding the distributed force influencing the flexible substrate that minimises a cost function, for example, the difference between the observed height of the thin film and some desired function. 

{\bf Novelty and contributions.}
The main contributions of our work are as follows:
\begin{enumerate}
\item 
We derive a thin-film flow model over a flexible substrate and introduce an externally forced substrate formulation that allows active control of the film through substrate deformation. The model couples the dynamics of the film and the flexible substrate, and includes the effects of gravity, disjoining pressure and external forcing on the substrate.

\item We prove that the system satisfies a global energy-dissipation law for a suitable energy functional. We also present an appropriate functional setting and a bounded weak formulation for our PDE system. The energy-dissipation law can be derived directly from the weak formulation, thereby providing an essential a~priori estimate on the PDE solution. 

\item  We derive the optimality conditions at the continuous level to control our coupled system. To achieve this, we formulate an optimal control problem consisting of the minimisation of a cost functional subject to the nonlinear evolution PDE for the thin-film flow controlled by the force acting on the substrate. First-order optimality conditions are given by a forward-backward, state-adjoint PDE system, closed with an optimality condition for the control.

\item We propose a stable time discretisation for the optimality system, which is possible because of the global energy-dissipation structure of the model.
In particular, a careful treatment of the nonlinearity in the model (i.e., the disjoining-pressure term) is needed to balance computational efficiency and stability.
%requires a particular choice in the time discretisation scheme 
For this, we propose a first-order Implicit-Explicit (IMEX) time-stepping scheme, which decomposes the nonlinear term into a contractive (convex) and expansive (concave) component~\cite{miles2020thermomechanically,wu2018posteriori}. This linearly-implicit scheme can be combined with any spatial discretisation method. We have employed standard continuous finite elements. Finally, gradient descent is used as an iterative method to solve the resulting discrete PDE-constrained optimisation problem~\cite{de2015numerical,manzoni2021optimal}. 
%The gradient descent algorithm is the most popular method which involves successive computations of optimality equations, with control updates formed following each iteration~\cite{manzoni2021optimal,de2015numerical}.
%Numerical results illustrate that controlling the thin film, even during rupture, enables achieving a precise film profile. This control strategy accelerates the system's achievement of a steady state, reduces instabilities, stabilises dewetting processes, and meets the desired profile specifications.
\end{enumerate}

The structure of the rest of the paper is as follows. Section~\ref{section2} introduces the model for the flow of thin films over a flexible substrate and its weak form. Section \ref{sec:control} shows the optimal control problem, the first-order necessary optimality conditions, and the gradient descent scheme for the optimisation problem. Section~\ref{sec:num_result2} presents the time discretisation using the IMEX time-stepping schemes and the numerical results. Finally, Section~\ref{Discussion} summarises the findings of this work.

%\lipsum[2-3]

% The outline is not required, but we show an example here.
\section{Model of thin-film flow over a flexible topography} \label{section2}

In this section, we present %our model which consists of a PDE of thin-film flow coupled to a flexible topography PDE. 
a mathematical model comprising two coupled PDEs that describe the evolution of the thin film and flexible substrate.
Furthermore, we show the model's functional settings and weak form.
The derivation of thin-film flow models on rigid substrates is well known and is based on lubrication theory~\cite{benney1966long,tseluiko2013stability}. The inclusion of the disjoining pressure and a non-inclined flexible substrate is discussed in various forms in the literature~\cite{matar2007dynamics,lee2011dynamics,matar2007falling, matar2004rupture,craster2009dynamics,alexander2020stability}.

\subsection{State equation} 
\label{sec1.2}
\begin{figure}[!t]
    \centering
\includegraphics[width=0.8\textwidth]{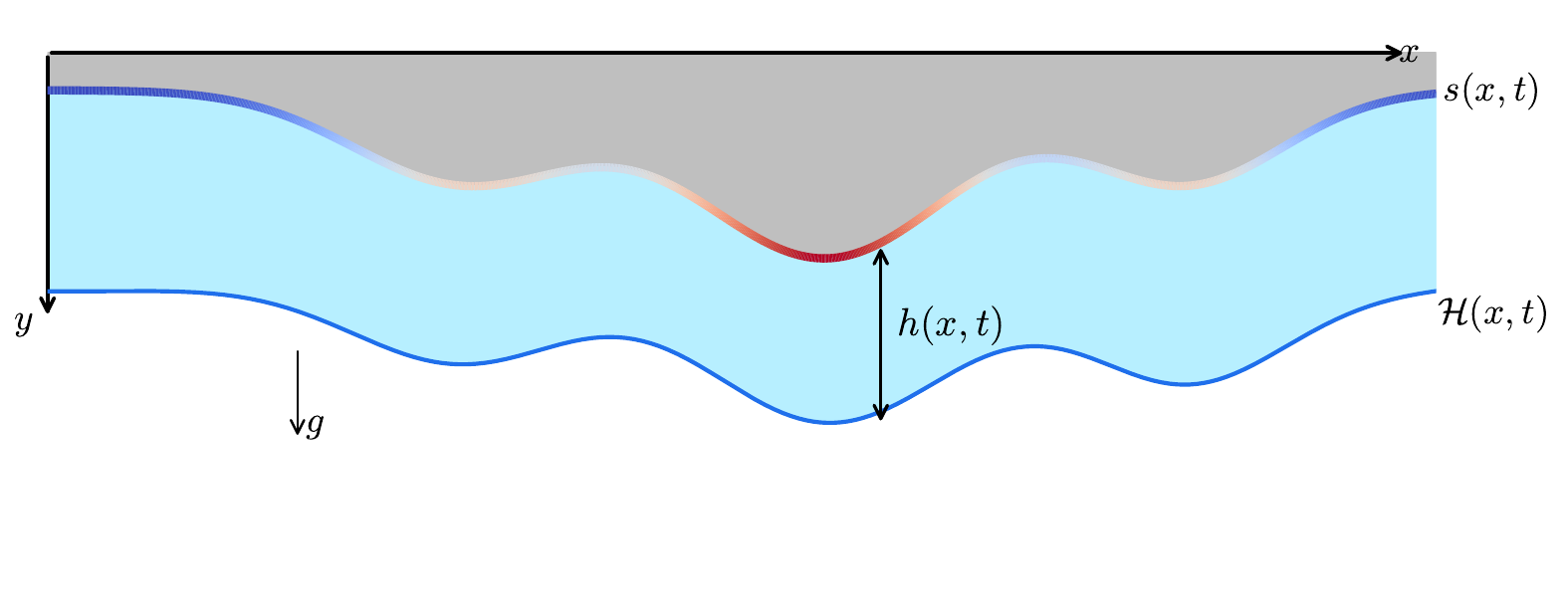}
  \caption{ Illustration of a two-dimensional fluid flow on a flexible substrate, where \( h(x, t) \) denotes the film thickness, \( s(x, t) \) is the shape of the topography, and \( \mathcal{H}(x, t) = h(x, t) + s(x, t) \) is the height of the one-dimensional free surface. The profile of the topography is controlled by a forcing function \(f(x,t)\), which is represented by the colour of the curve, where red indicates upward (positive) forcing and blue indicates downward (negative) forcing. \( g \) is the gravity vector.}
  \label{thin_film}
\end{figure}

We consider a 
thin-film flowing on the underside of a flexible substrate; an abstract sketch of this setup is illustrated in Figure~\ref{thin_film}. The coordinate system is defined by \((x,y) \), with $x$ being the horizontal coordinate 
and $y$ being the vertical coordinate. We also define $t$ as time. The force driving the flow is gravity $g$, which we take to act on the increasing \(y\) direction to model a film hanging beneath a substrate. 
This destabilised configuration generates rich dynamics in which our control can play a meaningful role.

The following system governs the evolution of a thin film over a flexible substrate, i.e.,  
\begin{equation}
\mathcal{H}(x,t)\coloneq h(x, t) + s(x, t),
\end{equation}  
where $\mathcal{H}(x,t)$ is the observed film height relative to flat horizontal reference, $h(x,t)\geq0$ is the film thickness, and $s(x, t)$ is the substrate topography height:
\begin{subequations}\label{Mixed}
\begin{align}
h_{t}&- \left(\tfrac{ 1}{3}h^3  \mu_{x}\right)_{x} = 0, \label{ht} \\
\mu &= \phi^{\prime}(h)-\frac{B_o}{C_a}\mathcal{H}-\frac{1}{C_a}\mathcal{H}_{xx}, \label{mu1} \\
s_{t} &- c^{2} s_{xx} -\frac{\gamma}{C_a} \left(h_{xx}+B_{o} h \right) = f(x,t). \label{flexible}
\end{align}
\end{subequations}
A derivation of the model equations is included in Appendix~\ref{app:derivation}. In equation~\eqref{ht}, we introduced  $\mu$ to reduce the fourth-order thin film PDE to a coupled system of second-order equations\footnote{Furthermore, when integrating $\mu$ using the finite element method, we find that the highest derivative is of the first order when we apply test functions.}. 
 
Equation~\eqref{mu1} includes an interface potential $\phi(h)$ capturing molecular interactions between liquid molecules and air \cite{duran2019instability,matar2004rupture}, which is commonly taken as 
\begin{equation}
\label{phi_h}
    \phi(h)=-\tfrac{A}{2}\,h^{-2}.
    \end{equation}
    This means there is a disjoining pressure $\Pi(h)=-\phi^{\prime}(h)=Ah^{-3}$ in the model to allow for the possibility of film rupture. Here, \(A\) is the dimensionless Hamaker constant.

%\delRone{ In \eqref{mu1}, the disjoining pressure $\Pi(h)=-\phi^{\prime}(h)$ is included to model film rupture, where $\phi(h)$ represents the interface potential capturing molecular interactions between liquid molecules and air \cite{duran2019instability} (see Section \ref{regul1} for examples of \(\phi(h)\)). }{} 
The second term on the right-hand side of \eqref{mu1} is associated with the action of gravity, and the last term reflects the effect of surface tension, where $C_a$ and $B_o$ denote the capillary and Bond numbers, respectively\footnote{The capillary and Bond numbers are dimensionless parameters common in fluid dynamics problems. They represent the importance of viscous forces or gravitational forces compared to surface tension forces on a fluid's surface, respectively. They are defined as \[C_a = \frac{\bar{\mu}U}{\sigma}, \quad B_o = \frac{\rho gh_0^2}{\sigma},\] where $\bar{\mu}$ and $\rho$ are the fluid's viscosity and density, respectively, $U$ is a typical velocity, $\sigma$ is the surface tension, $h_0$ is a typical length scale, and $g$ is gravity. See~\cite{kalogirou2020nonlinear} for typical values these parameters can take in real-life scenarios}.  
Furthermore, equation \eqref{flexible} represents a force balance on the substrate, containing a damping term $(s_{xx})$ that emerges from substrate tension, and a term $(h_{xx}+B_{o} h)$ that represents the pressure exerted on the substrate by the film. The parameters $c$ and $\gamma$ are dimensionless measures of the topography tension and film damping, respectively. The formulation of \eqref{flexible} also includes an external forcing term $f$, which enables direct control of the substrate dynamics.

To place our model~\eqref{Mixed} in context, we discuss how it relates with other models derived in the literature.
In the absence of disjoining pressure~\(\phi^{\prime}(h)=0\), and when the substrate is flat ($s\equiv 0$) or rigid, equations \eqref{ht}-\eqref{mu1} reduce to the well-established thin-film equation on a horizontal substrate, as used in \citep{hammond1983nonlinear,oron1997long} and \citep{jensen2004thin,tseluiko2013stability} for rigid substrates. Our model couples thin-film and flexible-substrate dynamics, and extends existing models in the literature (e.g.~\cite{matar2007dynamics, halpern1992fluid}) by additionally incorporating the effects of gravity, disjoining pressure and substrate forcing.
%in the following aspects. It is derived to additionally incorporate the effects of gravity, disjoining pressure and substrate forcing. 
The latter is particularly important in this work, as it allows the wall to be controlled actively. 

The system is solved over a time interval $t\in[0, T]$, where $T$ is the final time, and in a finite domain of arbitrary length $L$, i.e. $x \in \Omega \coloneq [0, L]$. The initial conditions are specified by setting
\begin{equation} \label{Ics}
     h(x, 0) = h_0(x), \quad s(x, 0) = s_0(x).
\end{equation}
Furthermore, homogeneous Neumann boundary conditions are applied at the domain boundaries
\begin{equation} \label{BC@}
    h_{x} = \mu_{x}=s_{x}=0,
    \qquad \text{at} \quad x=0 \text{ and } x=L.     
\end{equation}

We now clarify the free-energy dissipation in our model \eqref{Mixed}.  Understanding the energy dissipation mechanisms is important for ensuring the stability and physical relevance of the model. We begin by defining the total free-energy functional:
\begin{equation} \label{e_energy}
\mathcal{E}(h,s)=\int_\Omega\left[\phi(h)-\frac{B_o}{2 C_a} h^2+\frac{1}{2 C_a} h_x^2+\frac{c^2}{2 \gamma} s_x^2+ \frac{1}{C_a}(h_xs_x - B_o hs)\right] \mathrm{d}x,
\end{equation}
which consists of the molecular free energy \(\phi(h) - \frac{B_o}{2C_a} h^2\), surface free energy \(\frac{1}{2C_a} h_x^2\), elastic free energy \(\frac{c^2}{2 \gamma}s_x^2\), and fluid--substrate interaction energy \(\frac{1}{C_a}(h_xs_x - B_o hs\)). The latter is particularly important for the problem studied here, and describes how the combined variation of free surface and topography affects the energy -- more details are provided in Appendix \ref{energyApp}.

This free energy functional satisfies the following law (the full derivation is shown in Appendix \ref{energyApp}):
\begingroup
\small
\begin{equation} \label{energy}
\begin{aligned}
   \underbrace{ \frac{\mathrm{d}}{\mathrm{d} t} \mathcal{E}(h,s)}_{\text{Change in energy} }= \underbrace{ -\left(\frac{1}{3} \int_\Omega h^3 \mu_x^2 \,\mathrm{d}x + \frac{1}{\gamma} \int_\Omega s_t^2 \,\mathrm{d}x\right)}_{\text{Dissipation} } + \underbrace{\frac{1}{\gamma} \int_\Omega f s_t \,\mathrm{d}x}_{\text{External work}},
\end{aligned}
\end{equation}
\endgroup
where the terms on the right-hand side of \eqref{energy} are the dissipation in the fluid and substrate, and the work done by the external force on the substrate. Note that bracket in the dissipation term is non-negative (since \(h \geq 0\)), hence, in the absence of external work, the total free energy is non-increasing in time. 
We also note that the Bond number $B_o$ changes sign with $g$, but it affects only the potential-energy term in \eqref{e_energy} and not the dissipation.  
As a result, the system’s energy-decay property holds for either orientation of gravity.

It is well known that a free-energy dissipation law such as \eqref{energy} is critical in the study of well-posedness of thin-film flow models, see for example, \cite{bernis1990higher,klein2016optimal,bertozzi1994lubrication,binard2022well,gnann2018navier}. The well-posedness of system \eqref{Mixed} is outside of the scope of this work. We believe such a study could be applied to the natural weak formulation of the problem, which we present next.

\subsection{Weak form of the state equations}\label{Sec:weak}
We consider a mixed weak formulation of the state system \eqref{Mixed}, wherein $h$, $\mu$, and $s$ are considered as independent variables. 
We introduce  Hilbert spaces \(\hat{L}\) and \(\mathcal{V}\)  of functions dependent on time $t \in(0, T)$, possessing values in $L^2(\Omega)$ and $H^{1}(\Omega)$, respectively, and their norms:
\begin{subequations}
    \begin{alignat}{2}
        \hat{L} &\coloneq L^{2}\left(0, T ; L^{2}(\Omega)\right), & \quad\|v\|_{\hat{L}}^{2}\coloneq\int_{0}^{T}\|v(t)\|_{L^{2}(\Omega)}^{2} \,\mathrm{d} t, \\
        \mathcal{V} &\coloneq L^{2}\left(0, T ; H^{1}(\Omega)\right), & \quad\|v\|_{\mathcal{V}}^{2}\coloneq\int_{0}^{T}\|v(t)\|_{H^{1}(\Omega)}^{2} \,\mathrm{d} t. 
    \end{alignat}
\end{subequations}
The spaces $ \hat{L} , \mathcal{V}$ are fit for $f$ and $\mu$, respectively. A suitable space \(\mathcal{W}\) for $h$ and $s$ is defined by:
\begin{equation}
  \mathcal{W}\coloneq\left\{v \in \mathcal{V}: v_{t} \in \mathcal{V}^{*} \coloneq L^{2}\left(0, T ; H^{1}(\Omega)^{*}\right)\right\},  
\end{equation}
where $ H^{1}(\Omega)^{*}$ is the dual space of $H^{1}(\Omega)$, that is, the space of all continuous linear functionals on $H^{1}(\Omega)$. 
In addition, to strongly impose the initial conditions~\eqref{Ics} for $h$ and $s$, we establish
\begin{equation}
   \mathcal{W}_{v_{0}} \coloneq\left\{v \in \mathcal{W}: v(0)=v_0\right\},
\end{equation}
for any \(v_0 \in \hat{L}\).

The weak formulation of \eqref{Mixed} reads as follows (recall \(\mathcal{H}= h+s\)):
\par $ \text{Find}\ (h,\mu,s) \in \mathcal{W}_{h_{0}}  \times \mathcal{V} \times \mathcal{W}_{s_{0}} $, such that for %(almost every) 
\( t \in (0,T] \),
\begin{subequations}
\label{weak}
\begin{alignat}{2}
\left\langle h_t, v\right\rangle &+\frac{1}{3} (  h^3 \mu_x, v_x )  = 0,   & \forall\ v \in H^1(\Omega), \label{weak1} \\
(\mu, w)&-(\phi^{\prime}(h), w)+ \frac{B_o}{C_a}  (\mathcal{H}, w)-\frac{1}{C_{a}}\left(\mathcal{H}_x, w_x\right)=0, & \forall\ w \in H^1(\Omega), \label{weak2}  \\ 
\left\langle s_t, \psi\right\rangle &+c^2\left(s_x, \psi_x\right)+\frac{\gamma}{C_a} \left(\left(h_{x},\psi_x\right)- \left(B_{o} h,\psi \right) \right)=(f, \psi), & \quad \forall\ \psi \in H^1(\Omega)\label{weak3}. \end{alignat}
\end{subequations}
%\delRone{where $\langle\cdot, \cdot\rangle$ signifies the duality pairing between $H^{1}(\Omega)^*$ and $H^{1}(\Omega)$, and $(\cdot,\cdot)$ is the inner product in $L^2(\Omega)$.}{}

In \eqref{weak1}, the first term represents the duality pairing between
the time derivative and the test function,
$\langle h_t, v\rangle$.
Here, $\langle\cdot,\cdot\rangle$ denotes the pairing between
$H^{1}(\Omega)^{*}$ and $H^{1}(\Omega)$,
that is, the natural extension of the $L^{2}(\Omega)$ inner product
when the time derivative $h_t$ belongs to the dual space
$H^{1}(\Omega)^{*}$ rather than to $L^{2}(\Omega)$.
The notation $(\cdot,\cdot)$ is reserved for the standard
$L^{2}(\Omega)$ inner product~\cite[Sec.~2.1, Sec.~2.7]{de2015numerical,Troeltzsch2010}.

The system \eqref{weak} can also be written as an aggregated space--time weak formulation:
\begin{subequations} \label{weak_simple}
\begin{align} 
\text{Find}\ &(h,\mu,s) \in \mathcal{W}_{h_{0}}  \times \mathcal{V} \times \mathcal{W}_{s_{0}}: \notag \\
&\int_{0}^{T} \Big(\left\langle h_{t}, v\right\rangle +\left\langle s_{t}, \psi\right\rangle -(\phi^{\prime}(h), w) +(\frac{1}{3} h^3 \mu_x, v_x) \Big) \,\mathrm{d}t \\ 
&+\mathcal{A}\Big((h, s,\mu),(v, w,\psi)\Big)= \int^{T}_{0}(f,\psi) \,\mathrm{d}t, \nonumber  \qquad\qquad \forall\ (v, \omega, \psi) \in \mathcal{V} \times \mathcal{V} \times \mathcal{V}, 
\end{align}
where the bilinear form  $\mathcal{A}$ takes all linear terms:
\begin{equation}
\begin{aligned}
\mathcal{A}((h, s,\mu),(v, \omega,\psi))\coloneq\int_{0}^{T} \biggl((\mu, w)+ \frac{B_o}{C_a}  (\mathcal{H}, w)-\frac{1}{C_{a}}\left(\mathcal{H}_x, w_x\right) \\ +( c^{2} s_{x},  \psi_{x})+\frac{\gamma}{C_a} \left(\left(h_{x},\psi_x\right)- \left(B_{o} h,\psi \right) \right)\biggr) \,\mathrm{d}t.
\end{aligned}
\end{equation}
\end{subequations}

To derive the energy identity from the weak form, we select the test functions \(v=\mu(t),\ w=h_t(t)\), and \(\psi=s_t(t)\) in~\eqref{weak}, and multiplying~\eqref{weak2} by \(-1\) and~\eqref{weak3} by \(\frac{1}{\gamma}\), we obtain: 
\begin{subequations}\label{eq:weak_energy_tested_integral}
\begin{align}
\langle h_t,\mu\rangle
&= -\frac{1}{3}\int_\Omega h^3\,\mu_x^2\,\mathrm{d}x,
\\[2mm]
-\langle\mu,h_t\rangle
&\;+\;\frac{\mathrm{d}}{\mathrm{d}t}\int_\Omega \phi(h)\,\mathrm{d}x
\;-\;\frac{B_o}{C_a}\bigl(h+s,h_t\bigr)
\;+\;\frac{1}{C_a}\bigl(h_x+s_x,(h_t)_x\bigr)=0,
\\[2mm]
\;\frac{c^2}{\gamma}\bigl(s_x,(s_t)_x\bigr)
&
-\;\frac{B_o}{C_a}\bigl(h,s_t\bigr)
+\;\frac{1}{C_a}\bigl(h_x,(s_t)_x\bigr)= - \frac{1}{\gamma}\int_\Omega s_t^2\,\mathrm{d}x+ \frac{1}{\gamma}\bigl(f,s_t\bigr).
\end{align}
\end{subequations}
Adding these identities yields the energy-dissipation law~\eqref{energy}. By integrating~\eqref{energy} from \(0\) to \(t\), we then obtain the following  a~priori estimate for the PDE solution:
\begin{equation}
    \mathcal{E}(h(t),s(t)) + \int_0^t \int_\Omega \Big( \frac{h^3}{3} \mu_x^2 + \frac{1}{\gamma} s_t^2\Big) \mathrm{d}x \mathrm{d}t = \mathcal{E}(h_0,s_0)+\frac{1}{\gamma} \int_0^t \int_\Omega f s_t \mathrm{d}x \mathrm{d}t, \end{equation}
    where \(h_0, s_0 \in H^1(\Omega)\) are the initial conditions.

\subsection{Boundedness of weak form and disjoining pressure assumptions} \label{regul1}
We now discuss the boundedness of the weak formulation, in particular, suitable choices of the disjoining pressure \(\phi^{\prime}(h)\). First note that all the terms in the bilinear form \(\mathcal{A}(\cdot,\cdot)\) can be shown to be bounded within the indicated function spaces (by using Cauchy-Schwarz inequalities). Next, observe that the time derivative terms \( \langle h_t,v \rangle\) and \(\langle s_t, \psi \rangle \) are bounded because \( h_t,s_t \in H^1(\Omega)^*\) and \(v,\psi \in H^1(\Omega)\).  Furthermore, the term  \((\frac{1}{3} h^3 \mu_x, v_x)\) can be bounded as follows: 
\begin{equation}
        |( h^3 \mu_x, v_x)| \leq \|h||^3_{L^{\infty}(\Omega)}  \; \| \mu_x \|_{L^2(\Omega)} \, \|v_x\|_{L^2(\Omega)} \leq C \|h\|^3_{H^1(\Omega)} \| \mu\|_{H^1(\Omega)} \|v\|_{H^1(\Omega)}\,,
\end{equation}
for all~\(h,\mu,v\in H^1(\Omega)\), where in the last step we used a Sobolev embedding inequality (in~1D). 

Finally, to ensure that the term \( (\phi'(h),w) \) is bounded, 
we assume for simplicity that \( \phi\in C^{1}(\mathbb{R}) \) and satisfies the 
Lipschitz condition
\begin{equation}
    |\phi'(\xi)-\phi'(\zeta)| \le C_{\mathrm{L}}\,|\xi-\zeta|,
    \qquad \forall\,\xi,\zeta\in\mathbb{R}.
\end{equation}
Using this property, one obtains (see Appendix~\ref{proof_in} for the full proof)
\begin{equation}
    |(\phi'(h),w)| \le 
    C_{\mathrm{L}}\bigl(\|h\|_{L^{2}(\Omega)} + C \bigr)\,
    \|w\|_{L^{2}(\Omega)},
    \qquad \forall\, h,w \in L^{2}(\Omega)\supset H^{1}(\Omega),
\end{equation}
where the constant \(C\) can be written as 
\(C = C_{\mathrm{L}}|h_{\star}| + C^{*}\|1\|_{L^{2}(\Omega)}\), 
with \(h_{\star}\) an arbitrary reference value.

% \begin{delblockRone}{C4}
% A common physically-motivated choice for $\phi(h)$ in the literature is
% \begin{equation} \label{regul0}
%     \phi(h) = -\frac{A}{2}h^{-2}, \qquad
%     \phi'(h) = A h^{-3},
% \end{equation}
% where $A$ is the Hamaker constant.
% \end{delblockRone}

The choice of potential $\phi(h)=-\tfrac{A}{2}\,h^{-2}$ introduced in~\eqref{phi_h} does not satisfy the Lipschitz condition and is not bounded, which actually causes the model \eqref{Mixed} to break down when the film ruptures, i.e., \( h^{-1}(x) \to \infty \) for some \(x \in \Omega\). Regularisation of \(\phi(h)\) is commonly applied in the literature to ensure the solution can evolve beyond film rupture.  We adopt the approach proposed in \cite{GomZeeBOOK-CH2017} (see also \cite{miles2020thermomechanically}) which is as follows:
\begin{equation} \label{regul0}
    \phi(h) = \begin{cases} 
\dfrac{A}{2 \varepsilon^4} h^2 - \dfrac{A}{ \varepsilon^2} & \text{if } h < \varepsilon, \\[5pt] 
-\dfrac{A}{2} h^{-2} & \text{if } h \geq \varepsilon, 
\end{cases}
\end{equation}
where $\varepsilon > 0$ is a small regularisation parameter.  The idea is that the regularised $\phi(h)$ is convex quadratic for $h < \varepsilon$, while keeping its original form for $h \geq \varepsilon$, and ensuring that $\phi(h)$ has a bounded minimum at $h = 0$. This regularisation allows evolution beyond film rupture without the model breaking down.  
The adopted approach also ensures that $\phi(h)$ and $\phi^{\prime}(h)$ are continuous at all points, particularly at $h=\varepsilon$. 
The regularisation parameter $\varepsilon$ should be small to ensure the dynamics towards the moment of film rupture are accurate compared to the non-regularised case.  
% Following~\cite{miles2020thermomechanically}, we set $\varepsilon = 0.1$ in our computations, where the initial height of the film.%

\section{Optimal control problem} \label{sec:control}

We formulate an optimal control problem for the thin film flow dynamics over a flexible substrate \eqref{weak}. Our goal is to find a distributed force \( f(x,t) \) that influences the topography \( s(x,t) \) as in \eqref{flexible} so that the optimal solution approaches the target state \( \bar{h}(x) \). This optimal solution is either the height of the thin film \( h(x,t) \) or the total observed surface \( \mathcal{H}(x,t) = h(x,t) + s(x,t) \).

Formally, we seek a distributed force \( f(x,t) \)  that minimises the cost function \( \bar{\mathcal{J}}(h,s,f) \):
\begin{equation} 
\label{cost0}
\bar{\mathcal{J}}(h,s,f) \coloneq \frac{1}{2}\left\|(h+\beta s)(\cdot,T)-\Bar{h}(\cdot)\right\|_{L^{2}(\Omega)}^{2} + \frac{\alpha}{2} \int_{0}^{T}\|f(\cdot,t)\|_{L^{2}(\Omega)}^{2} \,\mathrm{d}t, 
\end{equation} 
where $\beta \in \{0,1\}$ and $\alpha >0$ is a regularisation parameter. The first term in equation \eqref{cost0} measures the difference between the current state and the target state $\bar{h}$, while the second term, determines a penalty for the control energy. In the unconstrained control setting of this paper, this control penalty is crucial for ensuring the problem's well-posedness \cite{herzog2010algorithms}. We formulate the optimal control problem as:
\begin{align}\label{ocp}
    \underset{(h,s,f)}{\min}\;\bar{\mathcal{J}}(h,s,f)\,,\quad\text{subject to \eqref{weak}}\,.
\end{align}

\subsection{First-order optimality conditions}

In this work, we follow an ``optimise-then-discretise'' approach that generally ensures computational robustness and mesh independence, effectively addressing numerical challenges in problems with low regularity multipliers and improving convergence behaviour \cite{de2015numerical}. 
 We begin by deriving first-order optimality conditions for \eqref{ocp}. For this, we define the Lagrangian
\begin{equation} \label{lagra}
\begin{aligned}
\mathcal{L}((h,\mu, s),f, (p,q,r)) &\coloneq \bar{\mathcal{J}}(h,s,f) + \int_{0}^{T}\left\langle h_{t}, p\right\rangle +\left\langle s_{t}, r\right\rangle \mathrm{d} t\\
+&\int_{0}^{T}(\frac{1}{3} h^3 \mu_x, p_x)-(f,r)-(\phi^{\prime}(h), q) \mathrm{d} t+ \mathcal{A}\biggl((h, \mu,s),(p, q,r)\biggl).
\end{aligned}
\end{equation}
%In equation \eqref{lagra}, $\mathcal{L}$ merges the cost function \eqref{cost} and the constraint \eqref{Mixed} into a singular function. Moreover, it is linear in terms of the Lagrange multipliers $(p(x,t), q(x,t),r(x,t)) $, serving as weights that balance the cost function and constraint within the optimisation process. 
%We aim to develop a set of equations that can define its solution accurately. Essential conditions for first-order optimality should be met by the solution. 
The optimality conditions arise from 
stationarity conditions of this Lagrangian with respect to state, adjoint, and control variables (see Chapter 3 in \cite{de2015numerical}  and \cite{albi2017mean} for further details). Consequently, the adjoint problem consists in finding \((p,q,r) \in \mathcal{W}^{h(T)} \times \mathcal{V} \times \mathcal{W}^{\beta s(T)}\) solving
\begin{align}
 & \hspace{-0.7cm}-\int_{0}^{T} \left(\left\langle p_{t}, v\right\rangle +\left\langle r_{t}, \psi\right\rangle +(\phi^{\prime \prime}(h)q, w)-\left(h^2 \mu_x p_x, v\right)+(\frac{1}{3} h^{3} p_{x}, w_{x} )  \right)\mathrm{d} t \notag \\ & \hspace{4cm}+\mathcal{A}\big((v, \omega,\psi),(p, q,r)\big)= \int^T_0(r,\psi) \,\mathrm{d}t, \label{weak_adj}
\end{align}
where $\mathcal{W}^{v_{T}}\coloneq\{v \in \mathcal{W}: v(T)=v_{T}\}$. 
The optimality conditions \eqref{weak_adj} can be expressed as a forward–backwards system. In other words, the optimal state variables, $(h^*,\mu^*,s^*)$, and the corresponding optimal adjoint variables, $(p^*,q^*,r^*)$, satisfy the following strong form:

\begin{subequations}\label{optsys}
  \begin{align}
    & \left\{
    \begin{aligned}
      & h^{*}_{t}-\frac{1}{3}\left(h^{*3} \mu^{*}_{x}\right)_{x}=0, \\
      & \mu^{*} =\phi^{\prime}(h^*)-\frac{B_o}{C_a}(\mathcal{H}^{*})-\frac{1}{C_a}\left(\mathcal{H}^{*}_{xx}\right), \\
      & s_t^*-c^2 s_{xx}^*-\frac{\gamma}{C_a} \left(h_{xx}^*+B_{o} h^* \right)=f^*,
    \end{aligned}
    \right.
    \label{state} \\[1em]
    & \left\{
    \begin{aligned}
      & -p^{*}_{t}+h^{*2} \mu^{*}_{x} p^{*}_{x}+\frac{1}{C_{a}} \left(q^{*}_{xx}+ \gamma r^{*}_{xx} \right) 
      +\frac{B_{o}}{C_{a}} \left(q^{*} - \gamma r^{*} \right)-\phi^{\prime\prime}(h^{*}) q^*=0, \\
      & \left(-\frac{1}{3} h^{*3} p^{*}_{x}\right)_{x}+q^{*}=0, \\
      & -r_t^*-c^2 r_{xx}^*+\frac{B_o}{C_a} q^*+\frac{1}{C_a} q_{xx}^*=0,
    \end{aligned}
    \right.
    \label{adjoint} \\[1em]
    & \left\{
    \begin{aligned}
      & h^{*}_{x}=\mu^{*}_{x}=s^*=p_x^*=q_x^*=r_x^*=0, && \quad \text{on} \ \Omega \times[0,T], \\
      & h^*=h_0, \quad s^*=s_0, && \quad \text{on} \ \Omega \times\{t=0\}, \\
      & p^{*}=r^{*}=\big( \bar{h}-(h+\beta s) \big), && \quad \text{on} \ \Omega \times\{t=T\}.
    \end{aligned}
    \right.
    \label{bc}
  \end{align}
\end{subequations}
closed by the optimality condition
\begin{equation} \label{control}
\alpha f^*-r^*=0.
\end{equation}

 \subsection{Reduced gradient approach}
We solve the optimisation problem \eqref{ocp} by means of gradient descent. For this, the system dynamics \eqref{Mixed} are used to define a nonlinear map from $f$ to $h+\beta s$, leading to a reduced objective functional $\mathcal{J}$ expressed solely in terms of $f$: 
\begin{equation} \label{reduced_cost}
    \mathcal{J}(f)=\bar{\mathcal{J}}(h(f), s(f), f)=  \frac{1}{2}\left\|(h+\beta s)(f; \cdot,T)-\Bar{h}(\cdot)\right\|_{L^{2}(\Omega)}^{2} +\frac{\alpha}{2} \int_{0}^{T}\|f(\cdot,t)\|_{L^{2}(\Omega)}^{2} \,\mathrm{d}t.
\end{equation}
Here, $(h+\beta s)(f;\cdot,t)$ represents the solution of the system \eqref{Mixed} for a given $f$. Expressed in terms of the reduced objective functional, the optimal control problem can be reformulated as follows: 
\begin{equation} \label{cost}
    \min_{f} \mathcal{J}(f), \quad \text{subject to \eqref{Mixed}.}
\end{equation}

It can be shown (see, e.g., Chapter 4 of \cite{de2015numerical}), that the gradient of the reduced objective can be recovered from the optimality condition \eqref{control}, that is,
\begin{equation}
    \nabla \mathcal{J}(f)\coloneq\alpha f-r\,,
\end{equation}
where the adjoint variable $r$ is to be computed from the forward-backward system \eqref{optsys} which depends on $f$ through the state equation. This gradient is the used in the iteration 
\begin{equation}
    f^{k+1}=f^{k}-\lambda^{k}\nabla \mathcal{J}(f^k), \quad k=1,2,\cdots,
\end{equation}
where the stepsize $\lambda^{k}$ is computed by backtracking. An elementary version of the pseudocode for this iterative loop is outlined in Figure \ref{algo}. 

\begin{figure} 
\begin{center}
\fbox{\parbox{0.7\linewidth}{    \begin{algorithmic}[1]
    \Require  $\text{tol}>0, k_{\max }, f^{0}, \lambda$, $k=0$;
      \While{$\left\|\nabla \mathcal{J}\left(f^k\right)\right\|>\text {tol and} k<k_{\max }$}
      \State  Obtain $\left(h^{k},\mu^{k}, s^{k} \right)$ with $f^{k}$ from \eqref{state}.
      \State  Obtain $\left(p^{k},q^{k}, r^{k}\right)$  with $\left(h^{k},\mu^{k}, s^{k} \right)$ from \eqref{adjoint}.
      \State  Evaluate the gradient $\nabla \mathcal{J}\left(f^k\right)$ from \eqref{control}.
      \State  Update $f^{k+1}=f^k-\lambda^{k} \nabla \mathcal{J}\left(f^k\right)$
      \While{ $\mathcal{J}\left(f^{k+1}\right) > \mathcal{J}\left(f^k\right)$} 
       \State   Obtain $\left(h^{k+1},\mu^{k+1}, s^{k+1} \right)$ with $f^{k+1}$
      \State   Evaluate the cost $ \mathcal{J}\left(f^{k+1}\right)$
      \If{$\mathcal{J}\left(f^{k+1}\right) > \mathcal{J}\left(f^k\right)$}
        \State $\lambda^m = 0.5 \lambda^k, m \geq k$  
        \State  Update $f^{k+1}=f^k-\lambda^{m} \nabla \mathcal{J}\left(f^k\right)$
        \EndIf 
     \EndWhile
       \State  $k:=k+1$ 
  \EndWhile
\end{algorithmic}
}}
\end{center}
\caption{Gradient descent method  with backtracking line search.}
\label{algo}
\end{figure}

  \section{Numerical approximation and tests}
   \label{sec:num_result2}
    \subsection{IMEX-FEM discretisation of the model} \label{regul}
For the computational treatment, the optimality conditions \eqref{optsys}-\eqref{control} derived above at a continuous level are discretised in space using the standard Finite Element Method (FEM) for variables $h, s, f$, and $\mu$, as well as the adjoint variables $p, q$, and $r$, subject to homogeneous Neumann boundary conditions. 
More details on the FEM formulation are provided in Appendix~\ref{FEM}.
For the time discretisation we need more careful treatment as the nonlinearities in equations \eqref{weak_simple} and \eqref{weak_adj} require a particular choice for the time-discretisation scheme to balance computational efficiency and stability \cite{wu2018posteriori}.
We explore first-order Implicit-Explicit (IMEX) time-stepping schemes, which decompose the nonlinear term $\phi$ in \eqref{regul0} as:
\begin{equation} 
\label{IMEX}
  \phi(h; v) = \phi_{+}(h; v) + \phi_{-}(h; v),  
\end{equation}
where $\phi_{+}$ and $\phi_{-}$ denote the contractive (convex) and expansive (concave) components, respectively. The IMEX scheme treats the convex part implicitly and the concave part explicitly. 
For the rupture term $\phi_+ = \tfrac{1}{2\varepsilon^4} h^2$, the derivative $\phi_+'(h) = \tfrac{1}{\varepsilon^4} h$ is linear.  
Treating it implicitly therefore introduces no nonlinear coupling and requires only a linear solve at each time step, retaining computational efficiency while ensuring stability as $h \to 0$.

The time interval $[0, T]$ is divided into $N$ smaller intervals $0=t_0<t_1<\dots<t_N=T$, with each sub-interval $\mathcal{I}_{k+1}=[t_k, t_{k+1}]$ having a uniform time step $\Delta t=t_{k+1}-t_k$ for $k=0,1,\dots,N-1$. 
The time discretisation for problem \eqref{weak_simple} is: 
\par Find $(h_{k+1}, \mu_{k+1}, s_{k+1}) \in \mathcal{V} \times  \mathcal{V}\times \mathcal{V}$ such that for all $(v, \omega, \psi) \in \mathcal{V} \times  \mathcal{V}\times \mathcal{V}$
\begin{equation} \label{000}
\begin{aligned}
 &\left(\frac{h^{k+1}-h^k}{\Delta t}, v\right) +\left(\frac{s^{k+1}-s^k}{\Delta t}, \psi\right) -\left(\phi^{\prime}_{+}(h^{k+1}) ; v\right)-\left(\phi^{\prime}_{-}(h^k) ; v\right) \\
 &\qquad\qquad +\left(\frac{1}{3} (h^k)^3 \mu_x^{k+1}, v_x\right)+\mathcal{A}((h^{k+1}, s^{k+1}, \mu^{k+1}),(v, \omega, \psi)) = (f^{k+1},\psi),
\end{aligned}
\end{equation}
for $k=0,1, \ldots, N-1$, where the initial condition is
$$
\left(h, v\right)=\left(h_0, v\right) \quad \text{and} \quad \left(s, \psi\right)=\left(s_0, \psi\right), \quad  \forall\ v,\psi \in L^2(\Omega) .
$$
It can be shown that this scheme dissipates the total energy \eqref{e_energy} using techniques shown in \cite{wu2018posteriori}.

For the adjoint problem \eqref{weak_adj}, we introduce an IMEX time-stepping scheme to discretise backwards in time:  
\par Find $(p_k, q_k, r_k) \in  \mathcal{V} \times  \mathcal{V}\times \mathcal{V}$ such that for all $(v, \omega, \psi) \in  \mathcal{V} \times  \mathcal{V}\times \mathcal{V}$
\begin{equation}
\begin{aligned}
-\left(\frac{p^{k+1}-p^k}{\Delta t}, v\right) & -\left(\frac{r^{k+1}-r^k}{\Delta t}, \psi\right) -\left(\phi^{\prime \prime}_{+}(h^{k+1})q^{k+1} , w\right)- \left(\phi^{\prime \prime}_{-}(h^k) q^{k} , w\right) \\
& -\left((h^k)^2 \mu_x^{k+1} p_x^{k+1}, v\right) +\frac{1}{3}\left( (h^k)^3 p_x^{k+1}, w_x \right)-(r^{k+1},\psi) \\
& +\mathcal{A}\left((p^k, q^k, r^k), (v, \omega, \psi)\right)= 0,
\end{aligned}
\end{equation}
for $k=0,1,2, \ldots, N-1$, where the terminal condition is:
\begin{equation}
\begin{aligned}
    \left(p^{N}, v\right)+\left(r^{N}, \psi\right) &- \Delta t \left(\phi^{\prime \prime}_{+}(h^{N}), q^{N}\right) + \frac{\Delta t }{3} \left((h^{N})^3 p_x^{N}, v\right) 
    \\ &+ \Delta t \mathcal{A}((p^{N}, q^{N}, r^{N}), (v, \omega, \psi)) = (\bar{h}^{N}-(h+\beta s)^{N}, v).
\end{aligned}
\end{equation}

 \subsection{Numerical tests} \label{numerics}
The results were generated using the forward-backward optimisation algorithm (Figure~\ref{algo}) with a tolerance $\text{tol}=10^{-4}$ and initial guess $f^0=0$. 
All simulations were carried out using a space-time descritisation with $n=250$ nodes, a time step $\Delta t=0.05$, and the following fixed parameters: $C_a=B_{o}=1$, $c=0.1$, $\gamma=0$, and the regularisation parameter $\alpha=10^{-6}$. For computational simplicity, we set $\gamma = 0$ in all simulations. This removes the pressure-coupling term from the substrate equation~\eqref{flexible}. The control \(f\) then only directly influences the topography~\(s\) of the wall via the substrate equation, which in turn affects the film height \(h\) via the thin-film system~\eqref{ht}-\eqref{mu1}. 

%The wall then responds only to the control $f(x,t)$ and its own elasticity, while still influencing the fluid through its deformation, which appears in the observed surface height $\mathcal{H} = h+s$.
We solve the system \eqref{optsys} on \([0, L]\). The initial conditions used in the numerical simulations are, unless otherwise specified,
\begin{equation} \label{IC}
    h_0(x) = 1 + h_A \cos\left(\frac{ m\pi x}{L}\right), \qquad
    s_0(x) = 0,
\end{equation}
where the initial amplitude $h_A$ will generally be chosen between $[0,0.5]$ (with $h_A=0$ corresponding to a flat state), and the frequency parameter will be set to $m=1$ unless otherwise specified.
We note that the respective flow with a flat substrate evolves to a nontrivial state whenever $L>2 \pi$, as predicted by linear stability analysis \cite{hammond1983nonlinear}. In addition to utilising the cost function and the norm of its gradient to validate our solution, we also employ the $L_{\infty}$-norm as a measure of accuracy; it is defined as $\max _{t \in[0, T]}\left\|(h^*(t)+\beta s^*(t))-\bar{h}\right\|_{\infty}$. 

For each example, three sets of figures are presented. The first set of figures show the free-surface and substrate evolutions in a symmetrically extended domain $[-L,L]$; this is done for illustrative purposes and better visualisation of the results. More detailed results, including the time evolution of the forcing $f(x,t)$, are demonstrated in the rest of the figures in the computational domain $[0,L]$. Table~\ref{tab:numerical_examples} summarises the key parameters and comparison between uncontrolled and controlled states for all numerical examples.
\begin{sidewaystable}
\centering
\caption{ Summary of numerical examples (shaded rows indicate rupture cases).}
\label{tab:numerical_examples}
\resizebox{\textwidth}{!}{%
\begin{tabular}{|>{\columncolor{white}}c|l|l|l|l|}
\hline
\multicolumn{5}{|c|}{\textbf{Common Parameters:} $C_a = B_o = 1$ (except Ex. 4: $B_o=0$), $c = 0.1$, $\gamma = 0$, $\alpha = 10^{-6}$, $n = 250$, $\Delta t = 0.05$, and $f^0=0$} \\
\hline
\hline
% \multicolumn{5}{|c|}{\textbf{Initial Conditions:} All examples follow initial conditions \eqref{IC}, except Ex.2: $s_0 = -0.25[\tanh(\frac{x+c_1}{d}) - \tanh(\frac{x-c_2}{d})]$, Ex.4: $h_0 = 1 + 10\exp(-(2x)^2)$, $s_0 = 0$, Ex.7: $h_0 = \bar{h}_\text{rupt}$ (initially ruptured)} \\
% \hline\hline
\textbf{Example} & 
\textbf{Key Parameters} & \textbf{Uncontrolled State} & \textbf{Controlled State} & \textbf{Figures} \\ \hline
\textbf{1} & 
% \begin{tabular}[c]{@{}l@{}}$h_0 = 1 + 0.5\cos(\pi x/L)$ \\ $s_0 = 0$\end{tabular} & 
\begin{tabular}[c]{@{}l@{}}$L = 3\pi$, $T = 5$, $h_A=0.5$ \\ $\beta = 1$, $\phi(h) = 0$\end{tabular} & 
\begin{tabular}[c]{@{}l@{}}Quasi-steady state reached \\ at $T = 900$\end{tabular} & 
\begin{tabular}[c]{@{}l@{}}Target state solution  from  Hammond \cite{hammond1983nonlinear} \\ achieved at $T = 5$ \end{tabular} & 
Figs. 3-5 \\ \hline

\textbf{2} & 
% \begin{tabular}[c]{@{}l@{}}$h_0 = 1$ (flat) \\ $s_0 = -0.25[\tanh(\frac{x+c_1}{d})$ \\ $\quad\quad - \tanh(\frac{x-c_2}{d})]$\end{tabular} & 
\begin{tabular}[c]{@{}l@{}}$L = 3\pi$, $T = 5$, $h_A=0$ \\ $s_0 = -0.25[\tanh(\frac{x+c_1}{d}) - \tanh(\frac{x-c_2}{d})]$ \\ $d = -0.2$ \\ $\beta = 1$, $\phi(h) = 0$,  $c_1 = -0.35L$ \\ $c_2 = 0.65L$\end{tabular} & 
\begin{tabular}[c]{@{}l@{}}Quasi-steady state reached  at $T = 900$\end{tabular} & 
\begin{tabular}[c]{@{}l@{}}Target state achieved  at $T = 5$\end{tabular} & 
Figs. 6-8 \\ \hline

\textbf{3} & 
% \begin{tabular}[c]{@{}l@{}}$h_0 = 1$ (flat) \\ $s_0 = 0$\end{tabular} & 
\begin{tabular}[c]{@{}l@{}}$L = 5$, $T = 1$, $h_A=0$ \\ $\beta = 1$, $\phi(h) = 0$\end{tabular} & 
\begin{tabular}[c]{@{}l@{}}Remains flat \\ (no evolution)\end{tabular} & 
\begin{tabular}[c]{@{}l@{}}Wave-like shape at $T = 1$\end{tabular} & 
Figs. 9-11 \\ \hline

\textbf{4} & 
% \begin{tabular}[c]{@{}l@{}}$h_0 = 1 + 10\exp(-(2x)^2)$ \\ $s_0 = 0$\end{tabular} & 
\begin{tabular}[c]{@{}l@{}}$L = 10$, $T = 5$ \\  $h_0 = 1 + 10\exp(-(2x)^2)$ \\ $\beta = 1$, $\phi(h) = 0$\end{tabular} & 
\begin{tabular}[c]{@{}l@{}}Non-flat steady state \\ (bump profile)\end{tabular} & 
\begin{tabular}[c]{@{}l@{}}Flattened to $\bar{h} = 1$  at $T = 5$\end{tabular} & 
Figs. 12-14 \\ \hline

\textbf{5} & 
% \begin{tabular}[c]{@{}l@{}}$h_0 = 1 + 0.5\cos(3\pi x/L)$ \\ $s_0 = 0$\end{tabular} & 
\begin{tabular}[c]{@{}l@{}}$L = \frac{15}{2}\pi$, $T = 10$ ,$h_A=0.5$ \\ $\beta = 1$, $\phi(h) = 0$\end{tabular} & 
\begin{tabular}[c]{@{}l@{}}Nonlinear instability \\ develops\end{tabular} & 
\begin{tabular}[c]{@{}l@{}}Maintains initial shape \end{tabular} & 
Figs. 15-17 \\ \hline\hline

\rowcolor{gray!25}
\textbf{6} & 
% \begin{tabular}[c]{@{}l@{}}$h_0 = 1 + 0.5\cos(\pi x/L)$ \\ $s_0 = 0$\end{tabular} & 
\begin{tabular}[c]{@{}l@{}}$L = 3\pi$, $T = 30$, $h_A=0.5$,  \\ $\beta = 0$, $\phi(h) \neq 0$, $\varepsilon = 0.1$, $A = 0.03$ \end{tabular} & 
\begin{tabular}[c]{@{}l@{}}Ruptured state reached \\ at $T = 550$\end{tabular} & 
\begin{tabular}[c]{@{}l@{}}Ruptured target state at $T = 30$\end{tabular} & 
Figs. 18-20 \\ 
\rowcolor{gray!25}
\hline

\rowcolor{gray!25}
\textbf{7} & 
% \begin{tabular}[c]{@{}l@{}}$h_0 = \bar{h}_\text{rupt}$ \\ (initially ruptured) \\ $s_0 = 0$\end{tabular} & 
\begin{tabular}[c]{@{}l@{}}$L = 3\pi$, $T = 30$  \\ $\beta = 0$, $\phi(h) \neq 0$, $\varepsilon = 0.1$, $A = 0.03$\end{tabular} & 
\begin{tabular}[c]{@{}l@{}}Remains ruptured\end{tabular} & 
\begin{tabular}[c]{@{}l@{}}Bubbles merge into  uniform surface at $T = 30$\end{tabular} & 
Figs. 21-22 \\ 
\rowcolor{gray!25}
\hline

\end{tabular}%
}
\end{sidewaystable}
\subsection{Non-rupture scenarios} \label{Norpture} 
In this section, we set \( \beta = 1 \) in \eqref{reduced_cost} because our goal is to control the observed surface, defined as \( \mathcal{H} = h + s \).

\subsubsection{Accelerated transition to desirable surface shape}
In the first two examples, where no disjoining pressure is present (i.e.\ $\phi'\equiv 0$), the uncontrolled dynamics are obtained by solving the forward problem~\eqref{Mixed} until temporal variations in the free surface become negligibly small. In Example~1 we consider $f=s=0$, so that \eqref{Mixed} reduces to the classical thin-film equation of Hammond~\cite{hammond1983nonlinear}, and the resulting quasi-steady profile is exported as the target state $\bar h$. In Example~2 we repeat the same procedure but over a prescribed, non-flat substrate $s(x)$, again integrating the uncontrolled system until the solution reaches a quasi-steady regime. These quasi-steady profiles provide the reference states against which the controlled evolution is compared.

\subsubsection*{Example~1}
We begin the presentation by comparing the time required for uncontrolled and controlled free surfaces to reach a quasi-steady state. For the uncontrolled case (i.e., $s=f=0$), the forward system \eqref{Mixed} has been resolved by Hammond~\citep{hammond1983nonlinear} for $L=3\pi$, and reached a quasi-steady state after $T=900$ time units. 
We aim to observe how fast the optimal state $\mathcal{H}^*$ approaches the target state $\bar{h}$ when control is applied.

\begin{figure}[!p]
  \centering
  \hspace{-1.5cm}
  \begin{subfigure}[b]{0.4\textwidth}
    \centering    \includegraphics[width=0.8\linewidth]{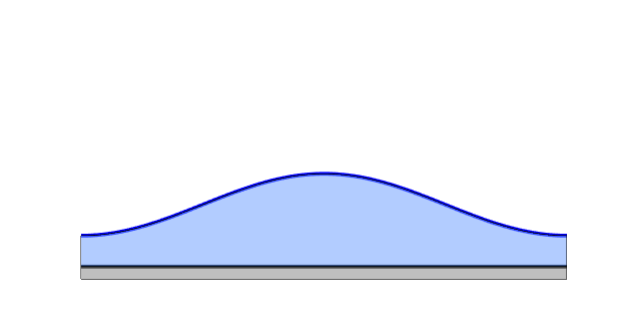} 
    \caption{$t=0$ }
  \end{subfigure}
  \hspace{-1.5cm}  
  \begin{subfigure}[b]{0.4\textwidth}
    \centering
\includegraphics[width=0.8\linewidth]{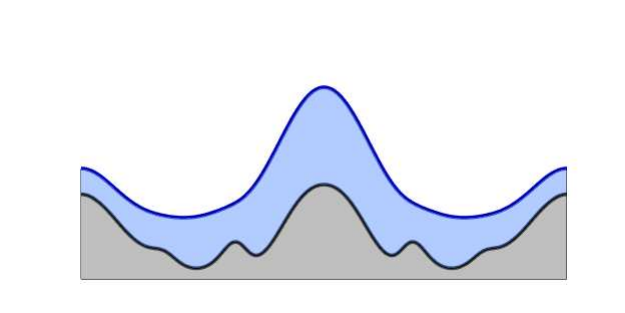}%{images/Hammond_T3.pdf}
\caption{$t=3$}
  \end{subfigure}
  \hspace{-1.5cm}  
  \begin{subfigure}[b]{0.4\textwidth}
    \centering
\includegraphics[width=0.8\linewidth]{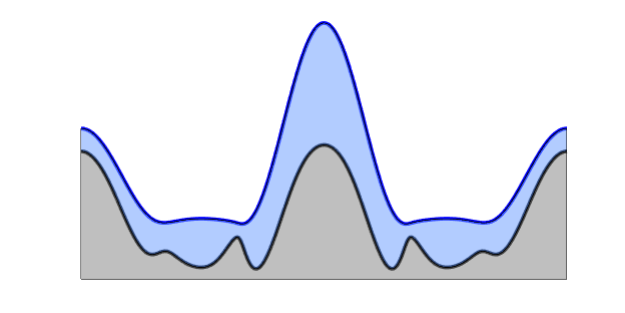}%{images/Hammond_T5.pdf}
\caption{$t=5$}
  \end{subfigure}
  \caption{Example~1. Snapshots of the optimal solution, controlled to a desired surface shape. The solution is extended symmetrically in domain $[-L,L]$.}
  \label{snapshots_Hammond}
\end{figure}

\begin{figure}[!p]
  \centering
   \hspace{-0.8cm}
  \begin{subfigure}[b]{0.4\textwidth}
    \centering
    \includegraphics[width=0.72\linewidth]{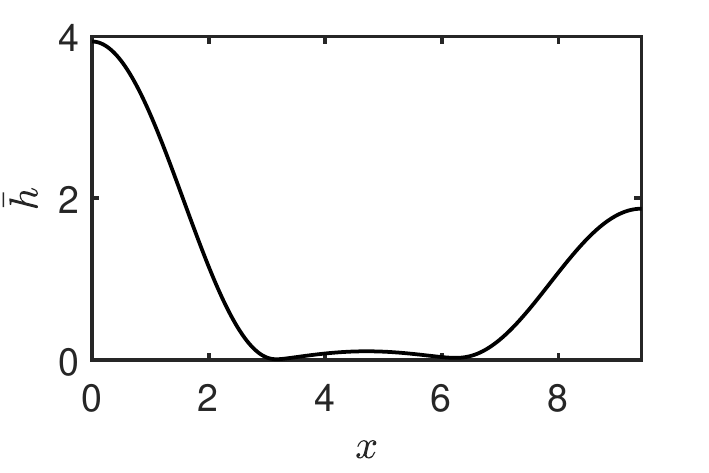}
    \caption{}
    \label{fig:sub1}
  \end{subfigure}
   \hspace{-1.5cm}
   \begin{subfigure}[b]{0.4\textwidth}
    \centering
    \includegraphics[width=0.72\linewidth]{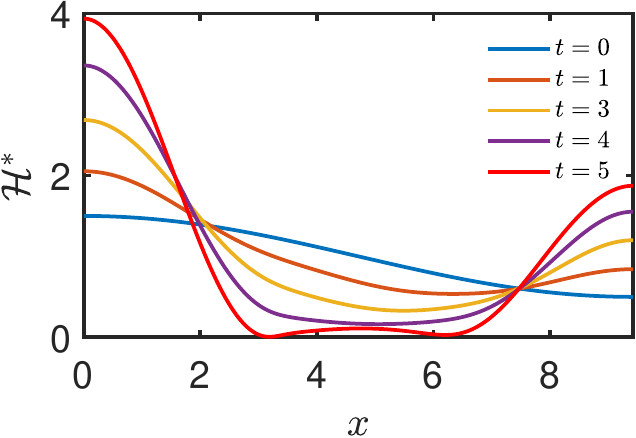}
    \caption{}
    \label{fig:sub3}
  \end{subfigure}
   \hspace{-1.3cm}
  \begin{subfigure}[b]{0.4\textwidth}
    \centering 
    \includegraphics[width=0.72\linewidth]{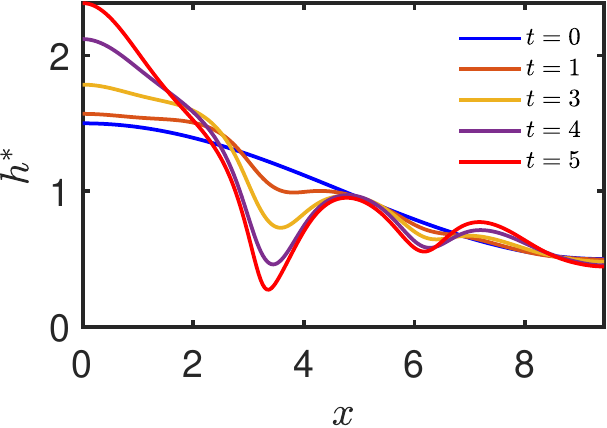}
    \caption{}
    \label{fig:sub_h}
  \end{subfigure}

    \hspace{-0.8cm}
  \begin{subfigure}[b]{0.4\textwidth}
    \centering
    \includegraphics[width=0.72\linewidth]{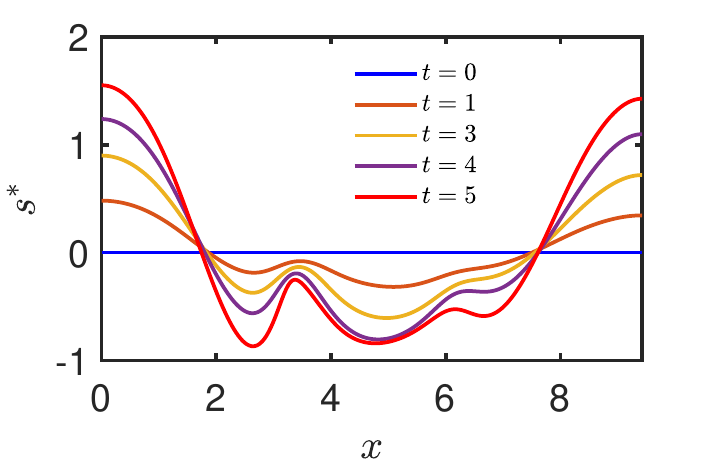}
    \caption{}
    \label{fig:sub4}
  \end{subfigure}
   \hspace{-1.5cm}
   \begin{subfigure}[b]{0.4\textwidth}
    \centering
    \includegraphics[width=0.72\linewidth]{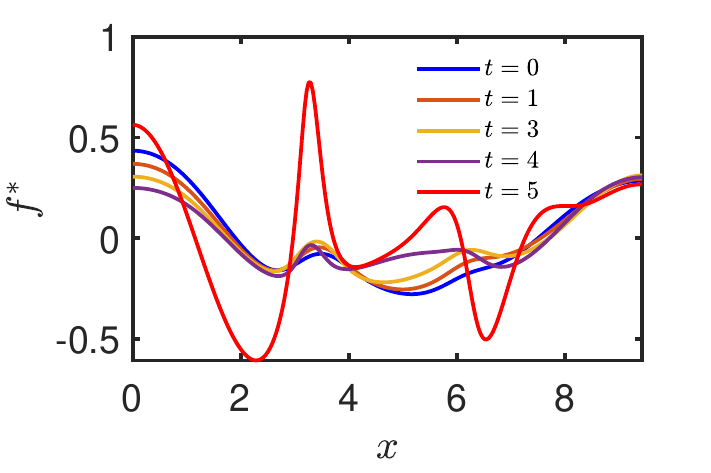}
    \caption{}
    \label{fig:hammond_f}
  \end{subfigure}
   \hspace{-1.3cm}
  \begin{subfigure}[b]{0.4\textwidth}
    \centering 
    \includegraphics[width=0.72\linewidth]{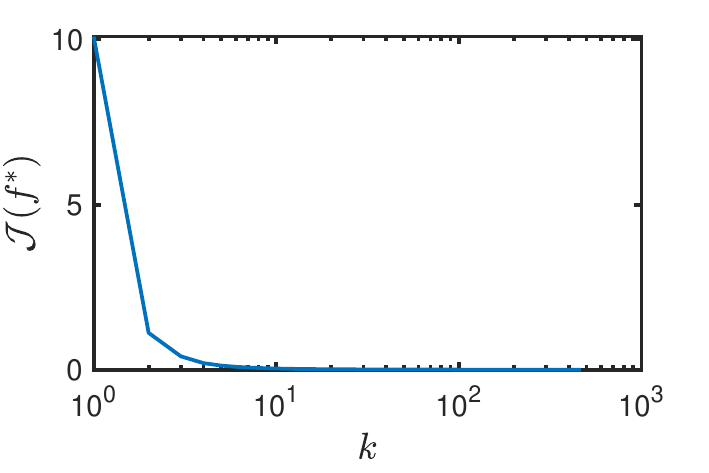}
    \caption{}
    \label{fig:sub6}
  \end{subfigure}
 \caption{%{\color{blue} Numerical results obtained using the parameter values $B_o=C_a=1$, $c=0.1$, $\gamma=0$, $\alpha=10^{-6}$, and $\beta=1$.}
 Example~1. (a) Target state $\bar{h}$. (b-e) Evolution of the optimal surface profile \( \mathcal{H}^* \), film height \( h^* \), topography \( s^* \), and control \( f^* \), respectively. (f) Cost function \(\mathcal{J}(f^*)\) against number of iterations \(k\).}
  \label{fig:speed_h}
\end{figure}

\begin{figure}[!p]
\centering
\hspace{-1.5cm}
\begin{subfigure}[b]{0.4\textwidth}  
    \centering
    \includegraphics[width=1.1\linewidth]{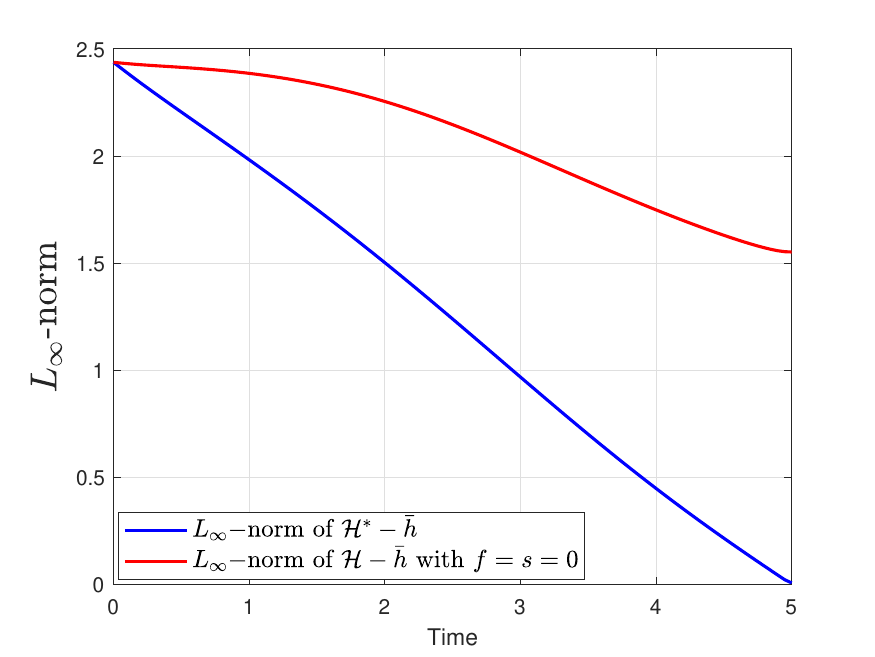}
    \caption{}
    \label{fig:sub5}
\end{subfigure}
\hspace{0.05\textwidth}  
\begin{subfigure}[b]{0.4\textwidth}  
\hspace{-1.0cm}
    \centering
    \includegraphics[width=1.1\linewidth]{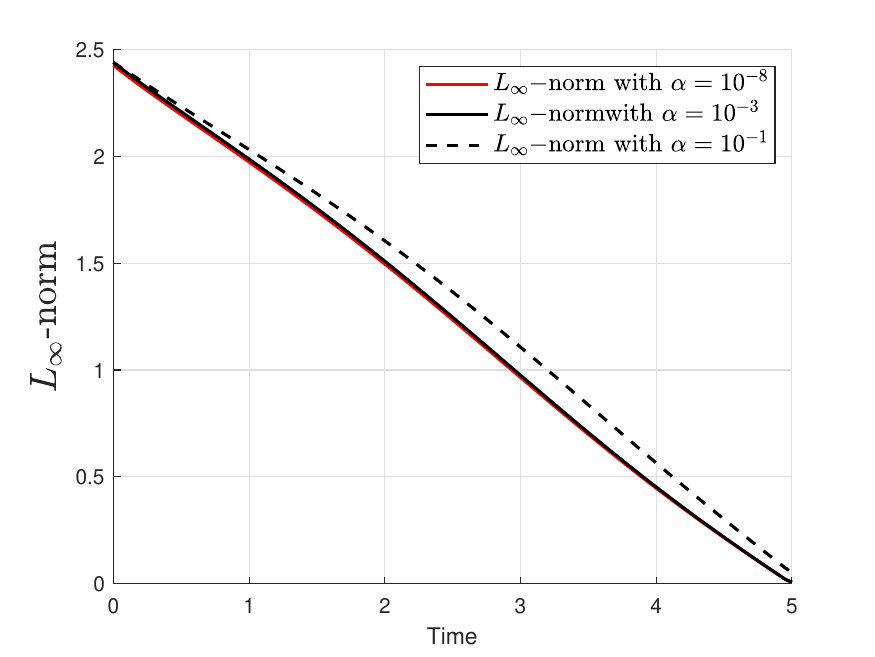}
    \caption{}
    \label{fig:hamm_alpha}
\end{subfigure}
\caption{Example~1. $L^{\infty}$-norms for (a) the controlled and uncontrolled solutions, and (b) various values of $\alpha$.}
\label{Fig3}
\end{figure}

We apply the forward-backwards optimisation method using the algorithm in
Figure~\ref{algo}, over the time horizon $[0,T]$ with $T=5$, using initial
conditions~\eqref{IC} with $h_A=0.5$. For illustrative purposes, snapshot plots spanning the domain $[-L, L]$ are presented in Figure \ref{snapshots_Hammond}.
The result of this computation is demonstrated in Figures~\ref{fig:speed_h} and \ref{Fig3}. 
The quasi-steady state is exported as our target state $\bar{h}$, shown in panel~\ref{fig:sub1}. 
The optimal free surface over topography, $\mathcal{H}^*$, shown in panel~\ref{fig:sub3}, successfully achieves this target state fast. 
The plot in panel~\ref{fig:sub_h} illustrates how the optimal free surface, $h^*$, varies over time. 
Panel~\ref{fig:sub4} shows the topography $s^*$ evolving from a flat state to an optimally adjusted profile with unexpected peaks and troughs, while panel~\ref{fig:hammond_f} shows the control parameter $f^*$ exhibits oscillatory behaviour to shape the optimal topography $s^*$. 
The cost function $\mathcal{J}(f^*)$ illustrated in panel~\ref{fig:sub6}, decreases to zero by the final time. 
Lastly, Figure~\ref{fig:sub5} shows the infinity norm of the controlled solution  $(h^*(t)+\beta s^*(t))-\bar{h}$ goes to zero much faster than the uncontrolled one, and Figure~\ref{fig:hamm_alpha} explores the impact of the regularisation parameter $\alpha$ on the solution, showing that smaller values of $\alpha$ than $10^{-3}$ lead to negligible changes to the optimal solution.

It is noteworthy to mention that, to ensure that $\mathcal{H}^*$ remains in a quasi-steady state,  the optimal solution  $\mathcal{H}^*$ and the topography $s^*$ are utilised as the initial conditions in the forward system \eqref{Mixed}. As expected, the free surface remains unchanged and the topography transitions towards a flat state.

\subsubsection*{Example~2} 
\begin{figure}[!p]
  \centering
  \hspace{-1.5cm}
  \begin{subfigure}[b]{0.4\textwidth}
    \centering    \includegraphics[width=0.8\linewidth]{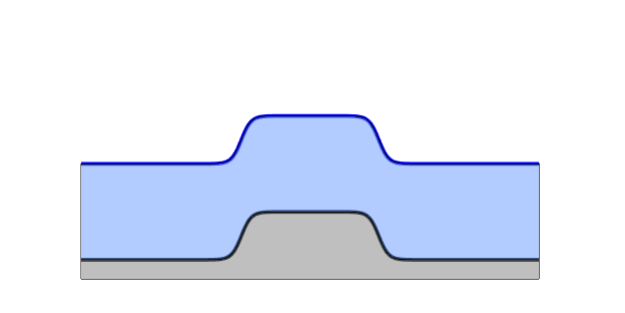}%{images/snap_0.pdf}
    \caption{$t=0$}
  \end{subfigure}
  \hspace{-1.5cm}  
  \begin{subfigure}[b]{0.4\textwidth}
    \centering  \includegraphics[width=0.8\linewidth]{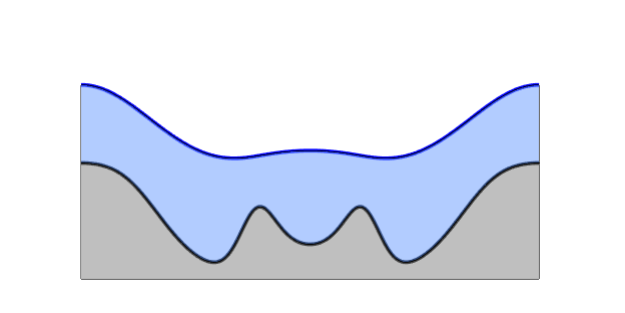}%{images/snap_3.pdf}
    \caption{$t=3$}
  \end{subfigure}
  \hspace{-1.5cm}  
  \begin{subfigure}[b]{0.4\textwidth}
    \centering
\includegraphics[width=0.8\linewidth]{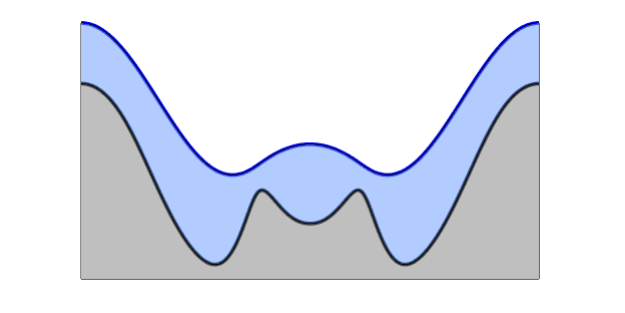}%{images/snap_5.pdf}
    \caption{$t=5$}
  \end{subfigure}
  \caption{Example~2. Snapshots of the optimal solution, controlled to a desired surface shape.}
  \label{snapshots_Hammond1}
\end{figure}

\begin{figure}[!p]
  \centering
   \hspace{-1.0cm}
  \begin{subfigure}[b]{0.4\textwidth}
    \centering
    \includegraphics[width=0.8\linewidth]{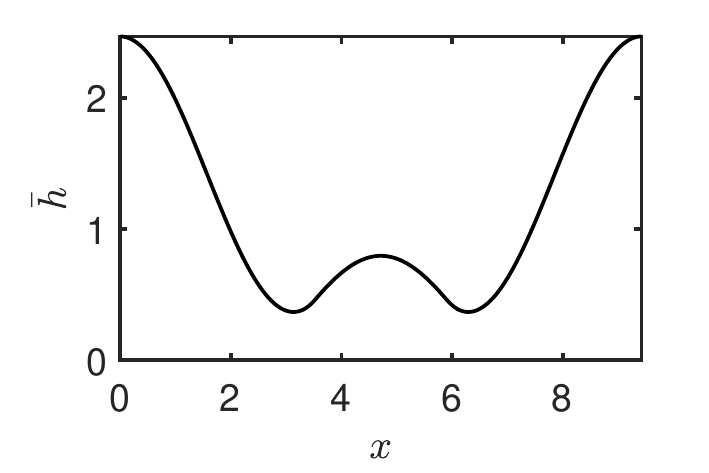}
    \caption{}
    \label{fig1:sub1}
  \end{subfigure}
   \hspace{-1.2cm}
   \begin{subfigure}[b]{0.4\textwidth}
    \centering
    \includegraphics[width=0.8\linewidth]{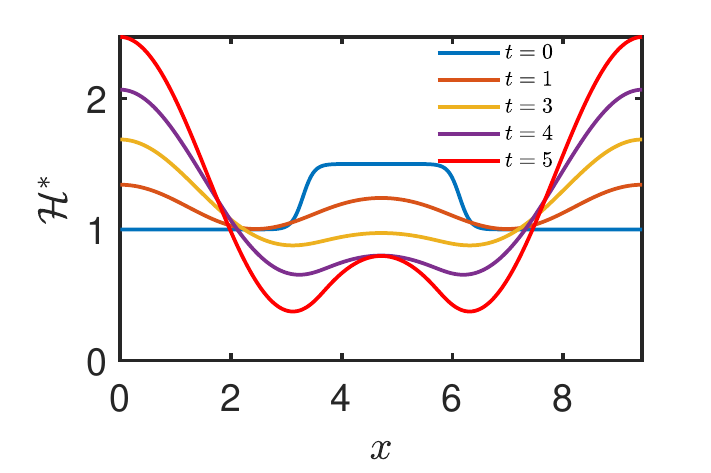}
    \caption{}
    \label{fig1:sub3}
  \end{subfigure}
   \hspace{-1.2cm}
  \begin{subfigure}[b]{0.4\textwidth}
    \centering 
    \includegraphics[width=0.8\linewidth]{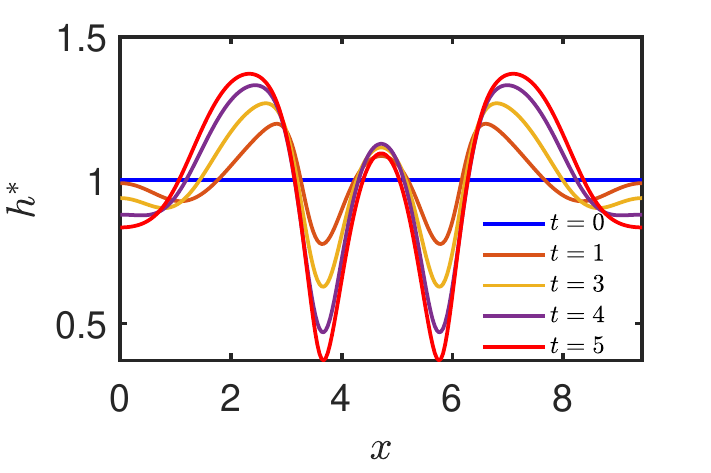}
    \caption{}
    \label{fig1:sub_h}
  \end{subfigure}
  
    \hspace{-1.0cm}
  \begin{subfigure}[b]{0.4\textwidth}
    \centering
    \includegraphics[width=0.8\linewidth]{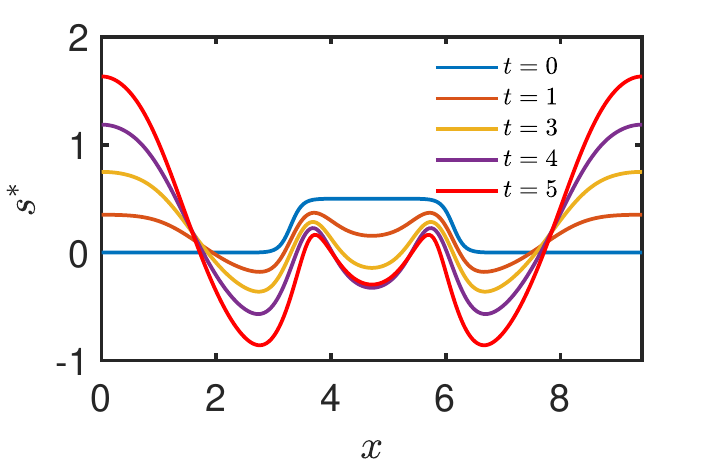}
    \caption{}
    \label{fig1:sub14}
  \end{subfigure}
  \hspace{-1.2cm}
  \begin{subfigure}[b]{0.4\textwidth}
    \centering
    \includegraphics[width=0.8\linewidth]{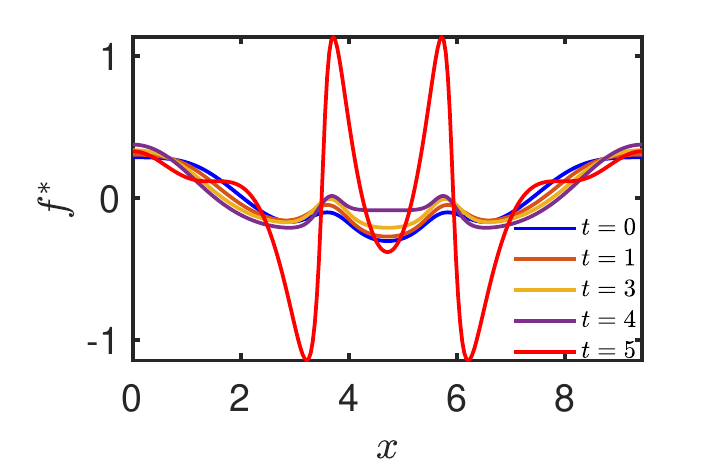}
    \caption{}
    \label{fig1:sub4}
  \end{subfigure}
\caption{Example~2. (a) Target state $\bar{h}$. (b-e) Evolution of the optimal surface profile \( \mathcal{H}^* \), film height \( h^* \), topography \( s^* \), and control \( f^* \), respectively.}
\label{Given_topo}
\end{figure}

\begin{figure}[!p]
\hspace{-2.5cm}
\begin{subfigure}[b]{0.4\textwidth}
    \centering
    \includegraphics[width=0.8\linewidth]{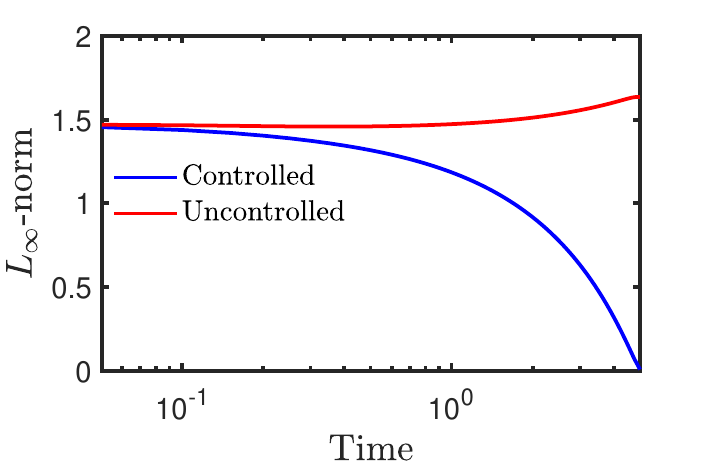}
    \caption{}
    \label{fig1:sub2}
  \end{subfigure}
\centering
  \begin{subfigure}[b]{0.4\textwidth}
    \centering 
    \includegraphics[width=0.8\linewidth]{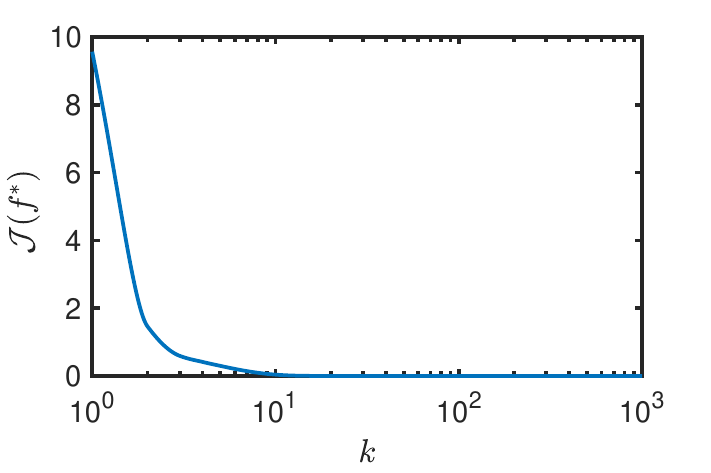}
    \caption{}
    \label{fig1:sub5}
  \end{subfigure}
\caption{Example~2. (a) \(L^\infty\)-norm for the controlled and uncontrolled solutions. (b) Cost function \(\mathcal{J}(f^*)\) against number of iterations \(k\).}
\label{givenT_cost}
\end{figure}

We extend the first example by solving the forward system \eqref{Mixed} over a non-flat topography, contrasting with the flat substrate examined previously. The uncontrolled system is solved with initial conditions comprising a flat surface profile (i.e., using $h_A=0$) and a substrate profile \( s_0(x) = a \left( \tanh \left( \frac{x + c_{1}}{d} \right) - \tanh \left( \frac{x - c_{2}}{d} \right) \right) \), where \( a = -0.25 \), \( c_{1} = -0.35L \), \( c_{2} = 0.65L \), $L=3\pi$, and the topography steepness is \( d = -0.2 \). 
The aim remains the same as in the first example in achieving the target state much faster than the uncontrolled one. A quasi-steady state for the free surface is achieved at \(T=900\) and the optimal control solution is presented in Figures~\ref{snapshots_Hammond1} and \ref{Given_topo}.
Panel~\ref{fig1:sub1} shows the target state, and panels~\ref{fig1:sub3} and \ref{fig1:sub_h} illustrate the time evolution of \(\mathcal{H}^*\) and \(h^*\), respectively. The smoothing changes in amplitudes and phases of \(s^*\) are highlighted in panel~\ref{fig1:sub14}, and the oscillatory forces \(f^*\) that influence \(s^*\) are shown in panel~\ref{fig1:sub4}. 
Figure~\ref{givenT_cost} demonstrates that both the \( L^\infty \) norm and the cost function converge to zero, although the convergence is slower than the case with a flat initial topography. 

%\subsubsection{ Desirable Surface Shape on a Short Time}
\begin{figure}[!p]
  \centering
   \hspace{-1.5cm}
  \begin{subfigure}[b]{0.4\textwidth}
    \centering
\includegraphics[width=0.8\linewidth]{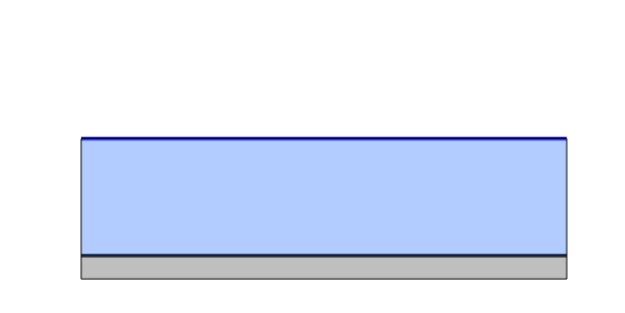}%
    %\caption{Profile of $\mathcal{H}^*-\bar{h}$}
    %\label{fig:sub2}
    \caption{$t=0$}
  \end{subfigure}
   \hspace{-1.5cm}
    \begin{subfigure}[b]{0.4\textwidth}
    \centering
\includegraphics[width=0.8\linewidth]{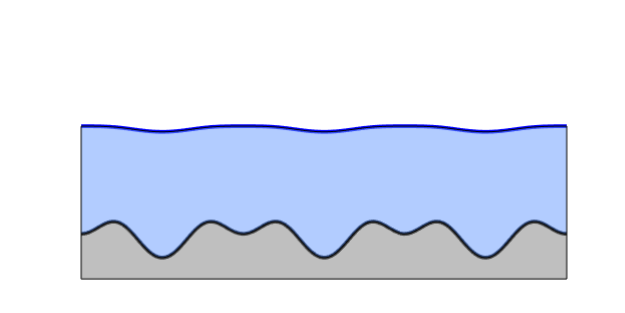}%{images/short_T025.pdf}
    %\caption{Evolution of optimal solution}
    %\label{fig:sub3}
    \caption{$t=0.5$}
  \end{subfigure}
   \hspace{-1.5cm}
  \begin{subfigure}[b]{0.4\textwidth}
    \centering
\includegraphics[width=0.8\linewidth]{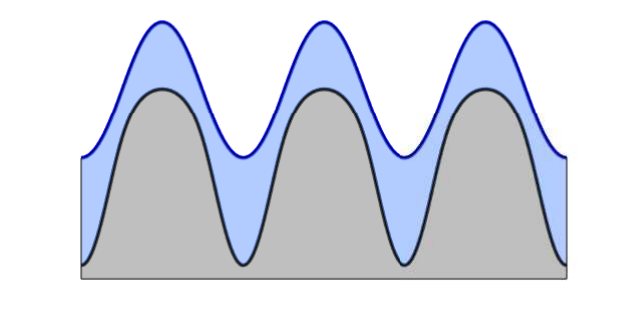}%{images/short_T05.pdf}
    %\caption{Evolution of optimal topography}
    %\label{fig:sub4}
    \caption{$t=1$}
  \end{subfigure}
  \caption{Example~3. Snapshots of the optimal solution, controlled to a desired surface shape. The solution is extended symmetrically in domain $[-L,L]$.}
  \label{snapshots_wave}
\end{figure}

\begin{figure}[!p]
\centering
 \hspace{-1.0cm}
 \begin{subfigure}[b]{0.4\textwidth}
    \centering
\includegraphics[width=0.8\linewidth]{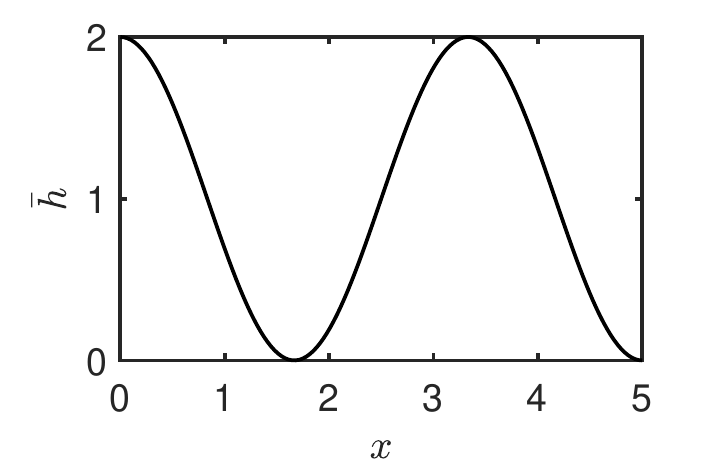}%{images/short_hbar.pdf}
    \caption{}
    \label{barwave}
  \end{subfigure}
\centering
 \hspace{-1.2cm}
  \begin{subfigure}[b]{0.4\textwidth}
    \centering
\includegraphics[width=0.8\linewidth]{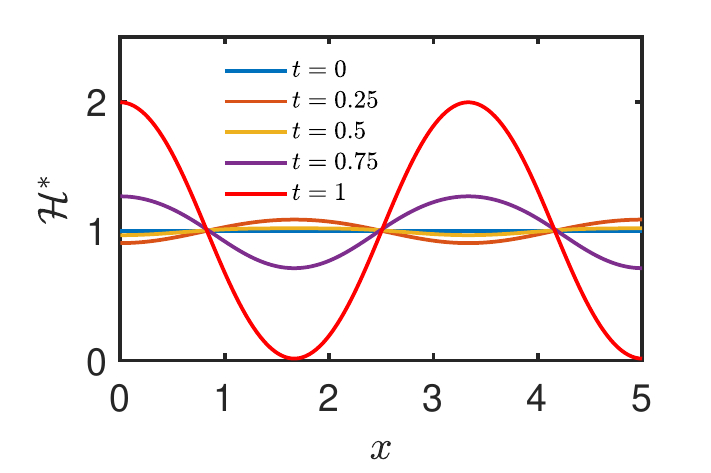}%{images/short_hbar.pdf}
    \caption{}
    \label{fig:wave_HS}
  \end{subfigure}
  \centering
 \hspace{-1.2cm}
  \begin{subfigure}[b]{0.4\textwidth}
    \centering
\includegraphics[width=0.8\linewidth]{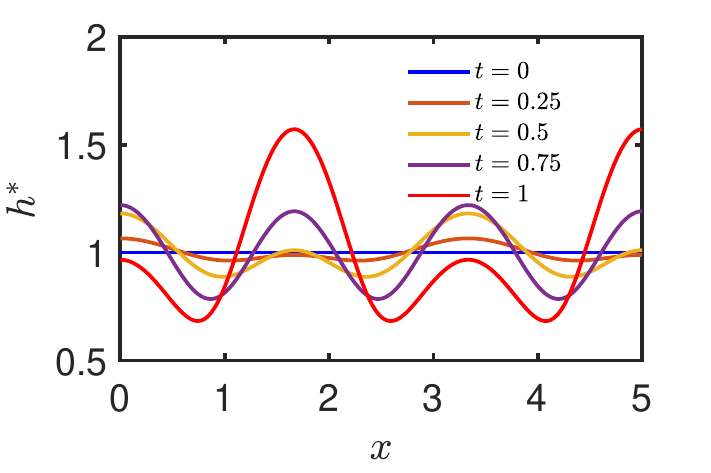}%{images/short_hbar.pdf}
    \caption{}
    \label{h_wave}
  \end{subfigure}
  \centering
   \hspace{-1.0cm}
  \begin{subfigure}[b]{0.4\textwidth}
    \centering
\includegraphics[width=0.8\linewidth]{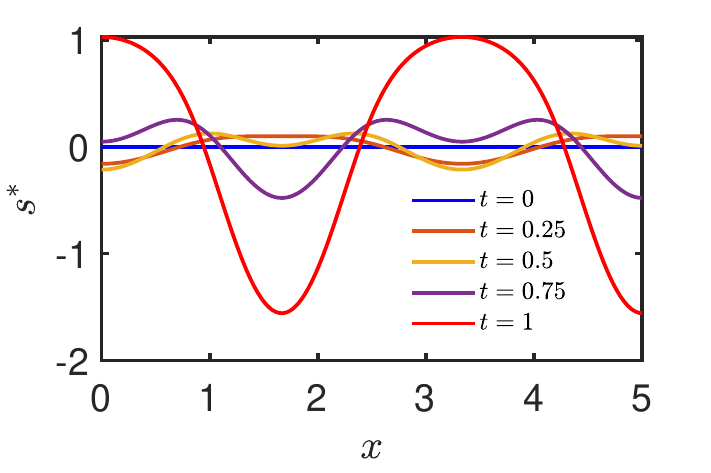}%{images/short_hbar.pdf}
    \caption{}
    \label{fig:wave_S}
  \end{subfigure}
  \centering
   \hspace{-1.2cm}
  \begin{subfigure}[b]{0.4\textwidth}
    \centering
\includegraphics[width=0.8\linewidth]{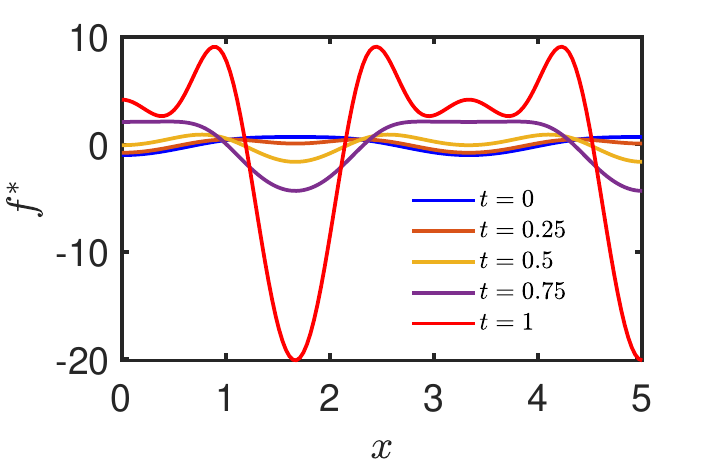}%{images/short_hbar.pdf}
    \caption{}
    \label{fig:wave_f}
  \end{subfigure}
  \caption{Example~3. (a) Target state $\bar{h}$. (b-e) Evolution of the optimal surface profile \( \mathcal{H}^* \), film height \( h^* \), topography \( s^* \), and control \( f^* \), respectively.}
  \label{fig:wave}
\end{figure}

\begin{figure}[!p]
    \centering
    \includegraphics[width=0.4\linewidth]{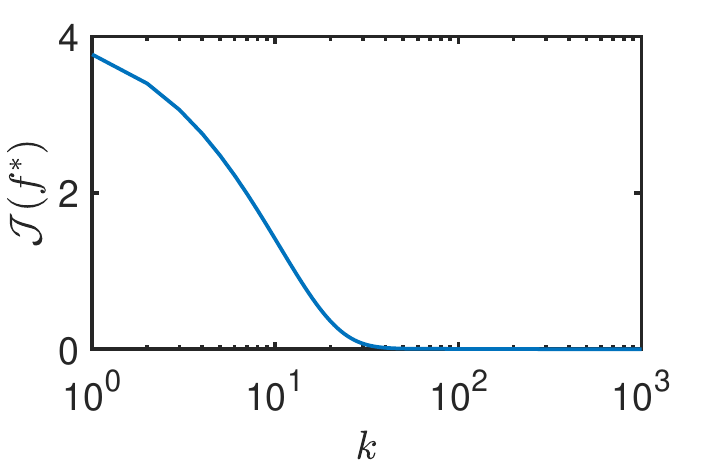}
    \caption{Example~3. Cost function \(\mathcal{J}(f^*)\) against number of iterations \(k\).}
    \label{fig:wave_cost}
\end{figure}

\subsubsection*{Example~3} The final example of this section considers a case where the free surface quickly transitions from flat to a desired wave-like shape in a short time $t=[0,1]$, as shown in Figures~\ref{snapshots_wave} and \ref{fig:wave}. 
Using a computational domain of length \(L=5\), we start with an initially flat surface by setting $h_A=0$ in \eqref{IC}. Without control, the surface remains flat, but with control, it transforms into the desired wave-like shape as shown in Figure \ref{barwave}. The optimal solution of \(\mathcal{H}^*\) and \(h^*\) are shown in panels~\ref{fig:wave_HS} and \ref{h_wave}, respectively. The topography \( s^* \) changes from flat to an oscillatory profile, enabling the desired wave-like target state. The control force \( f^* \) oscillates in sync with the topography changes, facilitating rapid adaptation, as seen in panel~\ref{fig:wave_f}. Figure~\ref{fig:wave_cost} demonstrates the cost function decreases towards zero, achieving the minimum cost, and the optimal control reaches the target surface.

\subsubsection{From non-flat to desired shapes}

\subsubsection*{Example~4} 
\begin{figure}[!p]
  \centering
   \hspace{-1.5cm}
  \begin{subfigure}[b]{0.4\textwidth}
    \centering
\includegraphics[width=0.8\linewidth]{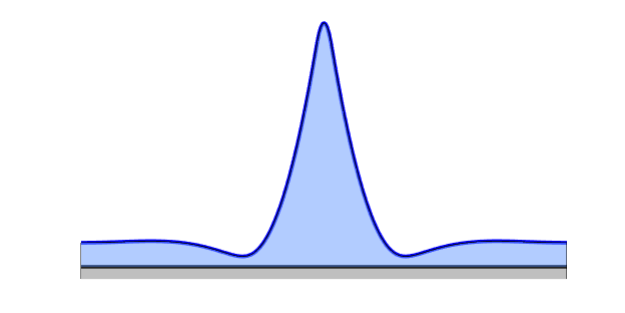}%{images/Jensen2_T0.pdf}
\caption{$t=0$}
  \end{subfigure}
   \hspace{-1.5cm}
  \begin{subfigure}[b]{0.4\textwidth}
    \centering
\includegraphics[width=0.8\linewidth]{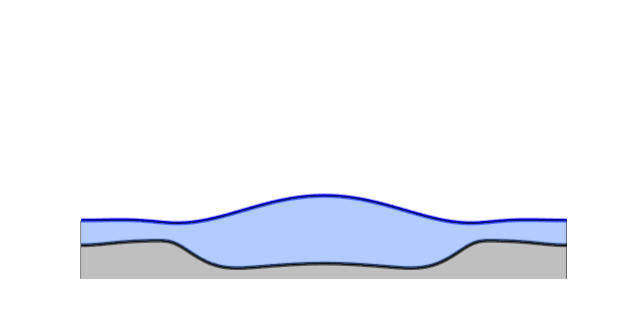}
\caption{$t=3$}
  \end{subfigure}
   \hspace{-1.5cm}
  \begin{subfigure}[b]{0.4\textwidth}
    \centering \includegraphics[width=0.8\linewidth]{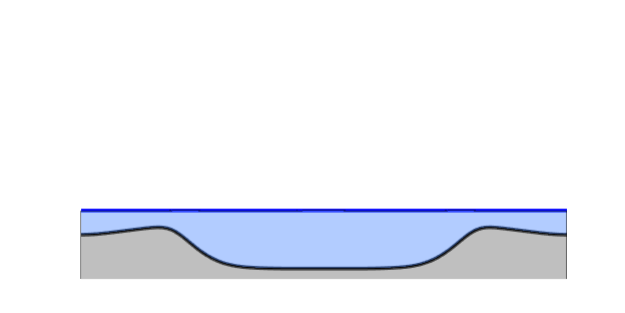}
\caption{$t=5$}
  \end{subfigure}
\caption{Example~4. Snapshots of the optimal solution, controlled to a uniform state. The solution is extended symmetrically in domain $[-L,L]$.}
\label{snapshots_Jensen}
\end{figure}

\begin{figure}[!p]
  \centering
   \hspace{-1.0cm}
   \begin{subfigure}[b]{0.4\textwidth}
    \centering    \includegraphics[width=0.7\linewidth]{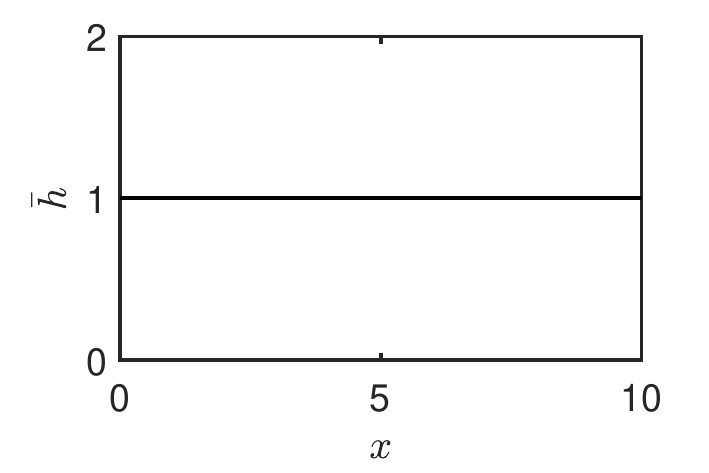}
    \caption{}
    \label{fig:jensen_hbar}
  \end{subfigure}
   \hspace{-1.2cm}
   \begin{subfigure}[b]{0.4\textwidth}
    \centering    \includegraphics[width=0.7\linewidth]{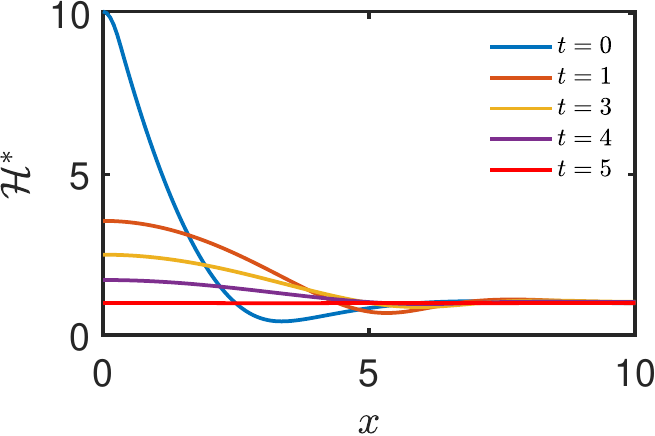}
    \caption{}
    \label{fig:jensen_HS}
  \end{subfigure}
   \hspace{-1.2cm}
  \begin{subfigure}[b]{0.4\textwidth}
    \centering 
    \includegraphics[width=0.7\linewidth]{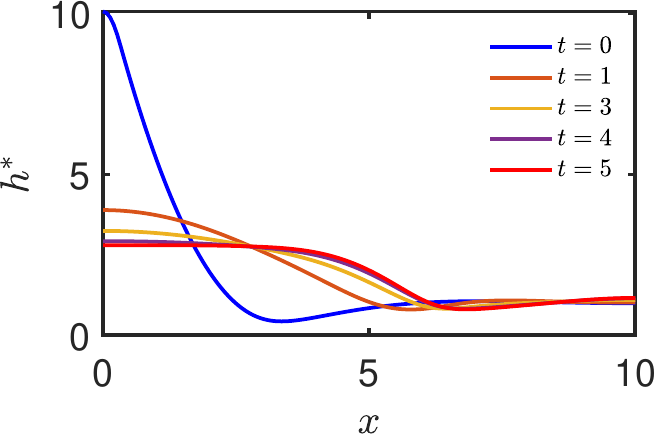}
    \caption{}
    \label{fig:sub_hjensen}
  \end{subfigure}
%\hspace{-1.5cm}
  \begin{subfigure}[b]{0.37\textwidth}
    \centering
    \includegraphics[width=0.8\linewidth]{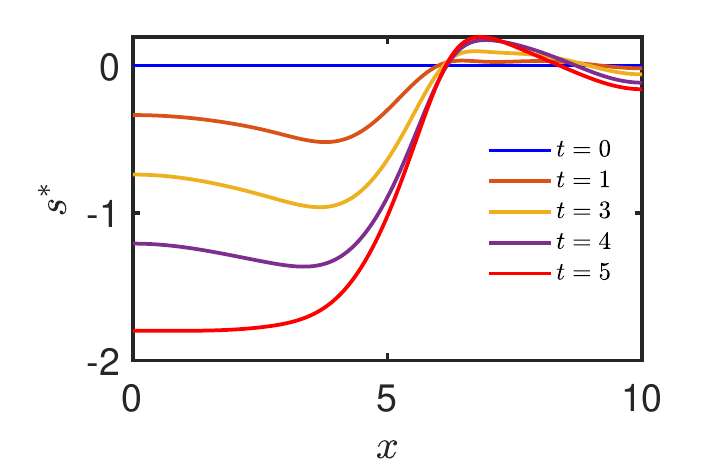}
    \caption{}
   \label{fig:jensen_s}
  \end{subfigure}
  \hspace{-1.2cm}
  \begin{subfigure}[b]{0.47\textwidth}
    \centering
    \includegraphics[width=0.7\linewidth]{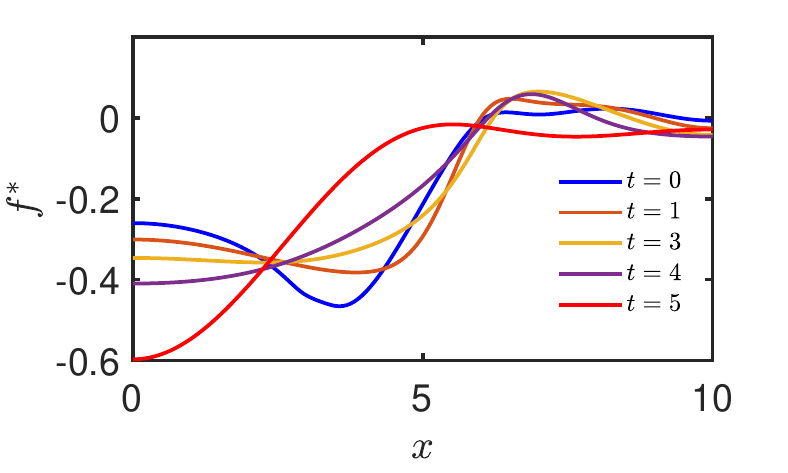}
    \caption{}
    \label{fig:jensen_f}
  \end{subfigure}
  \caption{Example~4. (a) Target state $\bar{h}$. (b-e) Evolution of the optimal surface profile \( \mathcal{H}^* \), film height \( h^* \), topography \( s^* \), and control \( f^* \), respectively.}
  \label{jensen results}
\end{figure}

The next example follows the work of \cite{jensen2004thin} where a related (uncontrolled) evolution equation for thin-film flow over topography was solved and a non-flat steady state was achieved. We initiate our control problem using one of these non-flat states as the initial condition by taking $h_0(x)=1+10 \exp \left(-(2x)^2\right)$, $s_0(x)=0$, and set \(L=10\), $B_o=0$. 
Our aim is to control the problem and flatten the surface (i.e., the target state is $\bar{h}=1$).
The time evolution of the various system variables is shown in Figures~\ref{snapshots_Jensen} and \ref{jensen results}.
Panel~\ref{fig:jensen_HS} shows that our optimal solution achieves a flat state, while panel~\ref{fig:sub_hjensen} displays deviations in \( h^* \) indicating it cannot be a flat alone. Panels~\ref{fig:jensen_s} and \ref{fig:jensen_f} reveal \( s^* \) and \( f^* \) adapting by \( T = 5 \) to counterbalance the variations in \( h^* \), aligning \( \mathcal{H}^* = h^* + s^* \) with the target state. 
Figure~\ref{fig:jensen_L_inf} shows the \( L^\infty \)-norm for the controlled and uncontrolled solutions over time, relative to the target state. The controlled solution moves towards zero, while the uncontrolled solution stabilises at a non-zero value.

\begin{figure}
   \begin{subfigure}[b]{0.4\textwidth}
    \centering
    \includegraphics[width=0.8\linewidth]{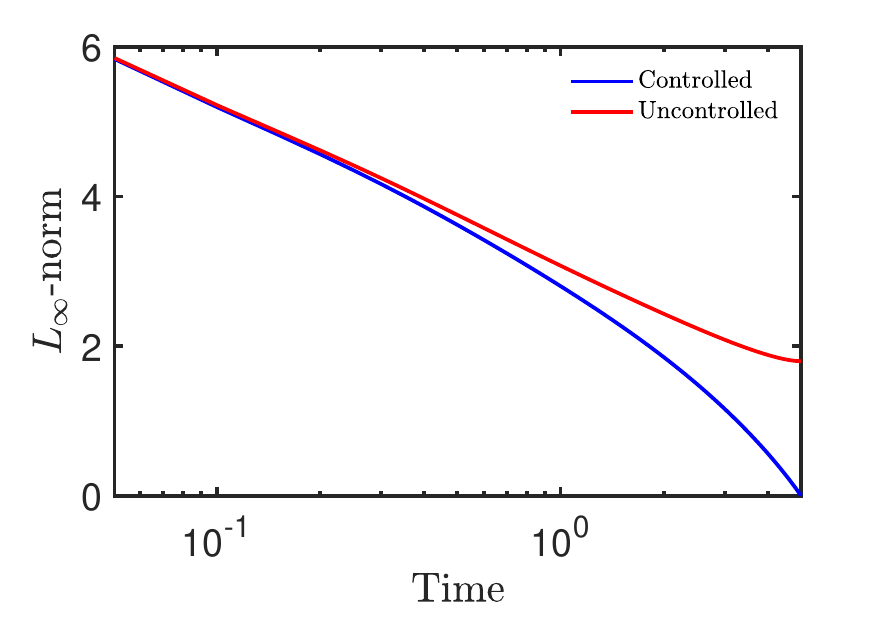}
    \caption{}
    \label{fig:jensen_L_inf}
  \end{subfigure}
  \begin{subfigure}[b]{0.4\textwidth}
    \centering    
    \includegraphics[width=0.8\linewidth]{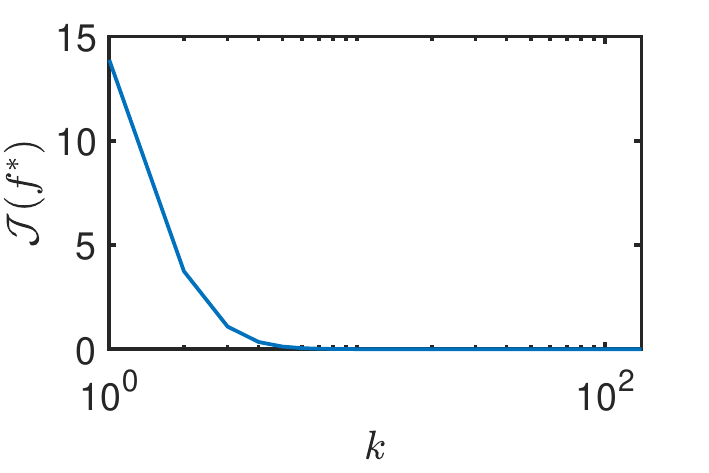}
 \caption{}
\label{fig:jensen_cost}
  \end{subfigure}
\caption{Example~4. (a) \(L^\infty\)-norm for controlled and uncontrolled solutions. (b) Cost function \(\mathcal{J}(f^*)\) against number of iterations \(k\).}
\label{fig:jensen_L_inf}
\end{figure}

\subsubsection*{Example~5} 
\begin{figure}[!p]
  \centering
  \hspace{-1.5cm}
  \begin{subfigure}[b]{0.4\textwidth}
    \centering
\includegraphics[width=0.8\linewidth]{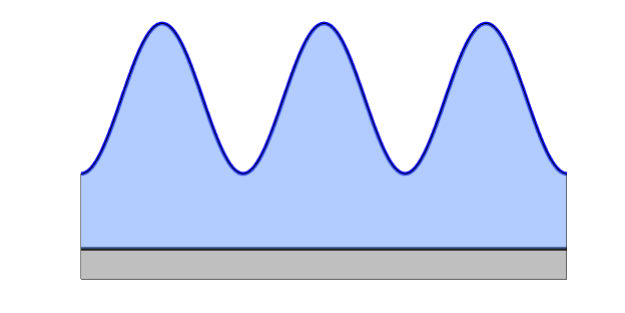}%
    \caption{$t=0$}
  \end{subfigure}
  \hspace{-1.5cm}
  \begin{subfigure}[b]{0.4\textwidth}
    \centering \includegraphics[width=0.8\linewidth]{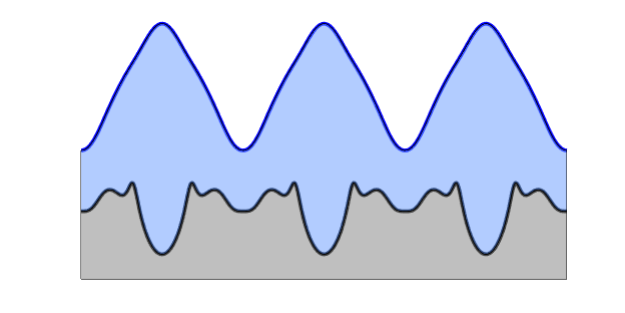}%{images/end_start_cT6.pdf}
    %\caption{$L_{\infty}$ Norm}
    %\label{fig:sub5}
    \caption{$t=5$}
  \end{subfigure}
  \hspace{-1.5cm}
  \begin{subfigure}[b]{0.4\textwidth}
    \centering    
    \includegraphics[width=0.8\linewidth]{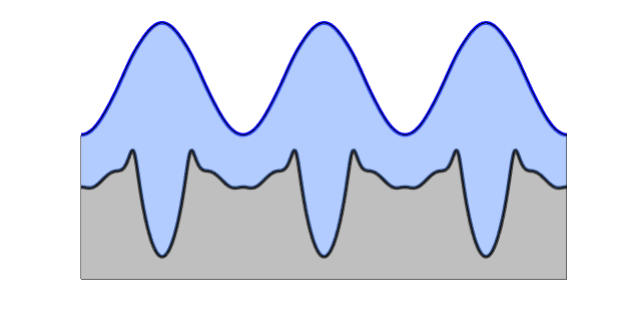}%{images/end_start_cT10.pdf}
    %\caption{Profile of $\bar{h}$}
   % \label{fig:sub1}
   \caption{$t=10$}
  \end{subfigure}
  \caption{Example~5. Snapshots of the optimal solution, controlled to remain at the initial surface shape. The solution is extended symmetrically in domain $[-L,L]$.}
  \label{snapshots_end}
\end{figure}

\begin{figure}[!p]
  \centering
  \hspace{-1.0cm}
  \begin{subfigure}[b]{0.4\textwidth}
   \centering \includegraphics[width=0.75\linewidth]{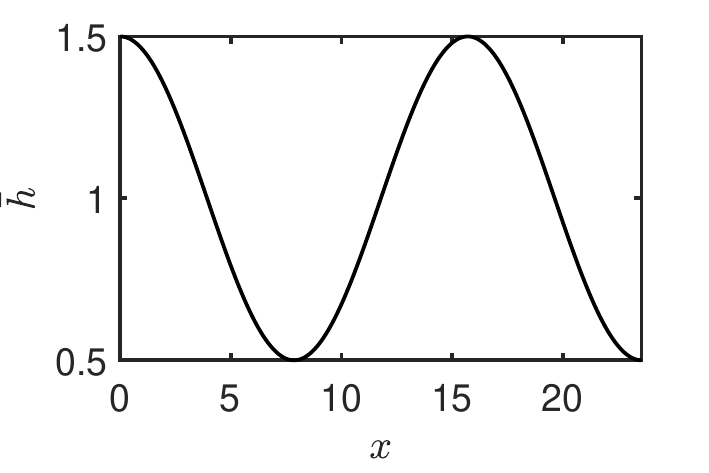}
    \caption{}
    \label{EX4_hbar}
  \end{subfigure}
   \hspace{-1.2cm}
   \begin{subfigure}[b]{0.4\textwidth}
    \centering   \includegraphics[width=0.72\linewidth]{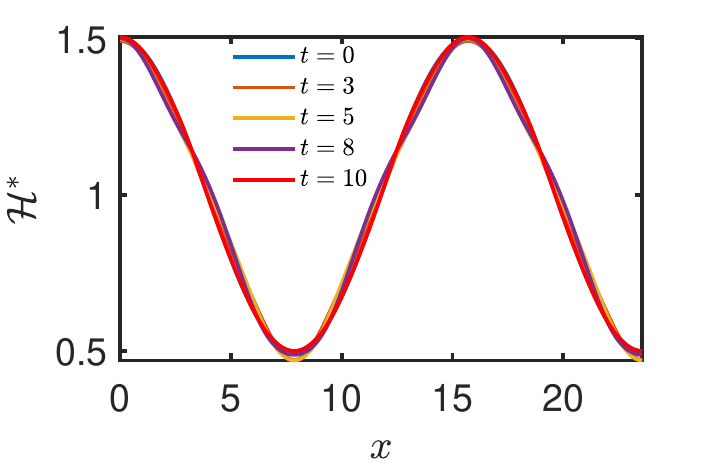}
    \caption{}
    \label{wave_H}
  \end{subfigure}
   \hspace{-1.2cm}
  \begin{subfigure}[b]{0.4\textwidth}
    \centering     \includegraphics[width=0.72\linewidth]{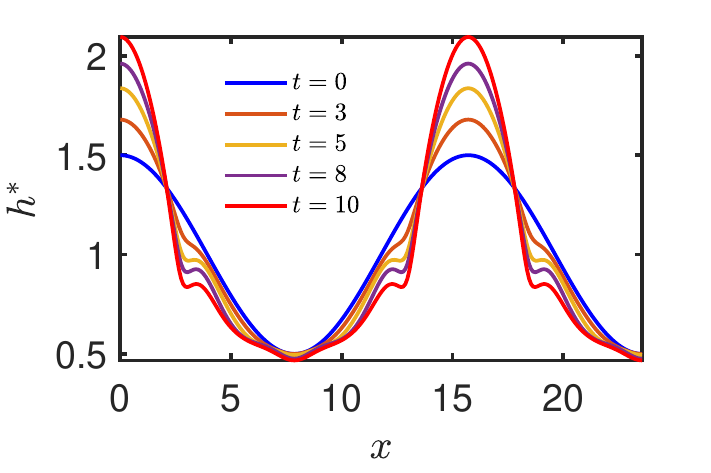}
    \caption{}
    \label{wave_h}
  \end{subfigure}
  
    \hspace{-1.0cm}
  \begin{subfigure}[b]{0.4\textwidth}
    \centering
\includegraphics[width=0.75\linewidth]{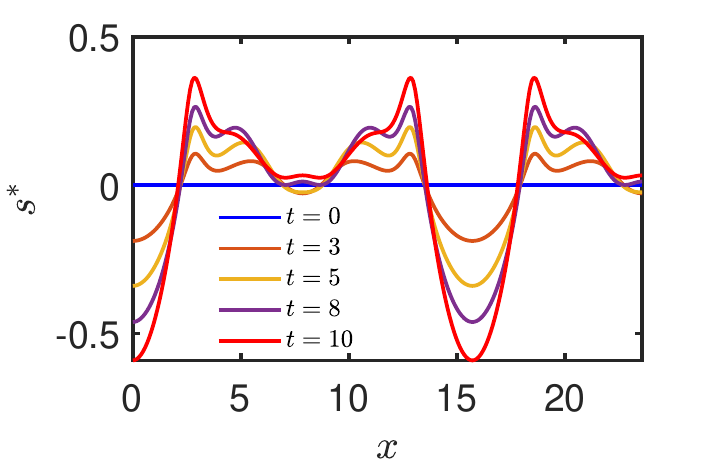}
    \caption{}
   \label{wave_s}
  \end{subfigure}
  \hspace{-1.2cm}
  \begin{subfigure}[b]{0.4\textwidth}
    \centering    
    \includegraphics[width=0.75\linewidth]{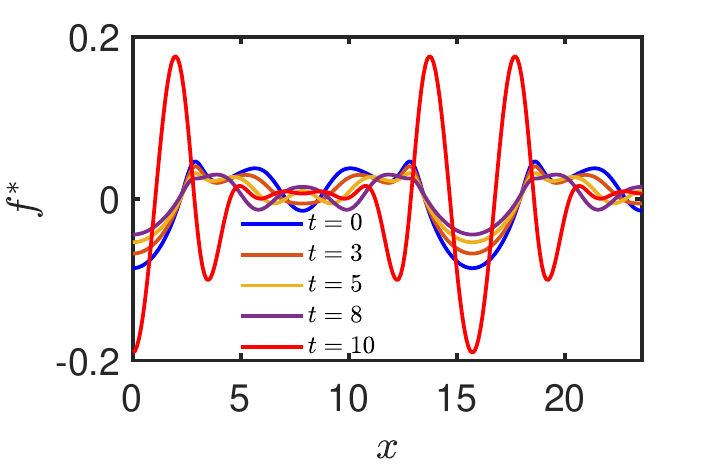}
    \caption{}
     \label{wave_f}
  \end{subfigure}
    \hspace{-1.2cm}
    \begin{subfigure}[b]{0.4\textwidth}
    \centering
\includegraphics[width=0.75\linewidth]{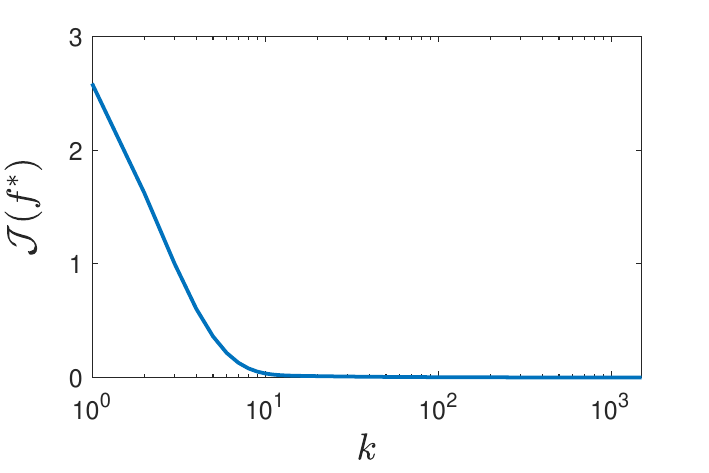}
\caption{}
\label{wave_cost}
  \end{subfigure}
\caption{Example~5. (a) Target state $\bar{h}$. (b-e) Evolution of the optimal surface profile \( \mathcal{H}^* \), film height \( h^* \), topography \( s^* \), and control \( f^* \), respectively. (f) Cost function \(\mathcal{J}(f^*)\) against number of iterations \(k\).}
\label{fig:wave2}
\end{figure}

\begin{figure}[!p]
    \centering
\includegraphics[width=0.4\linewidth]{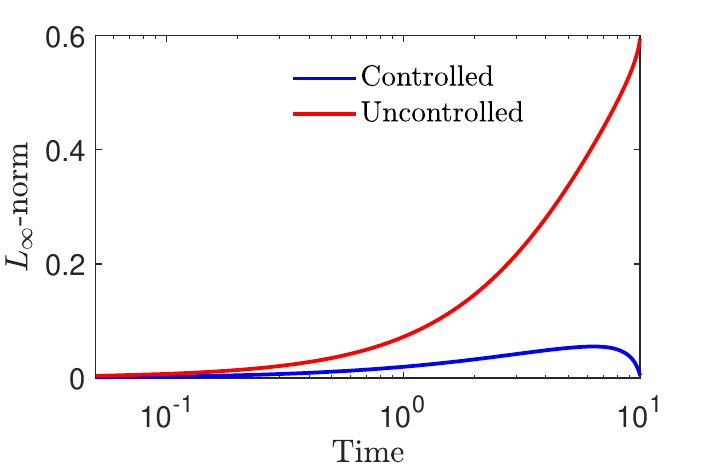}
   % \label{wave_norm_cost}
    \caption{Example~5. $L^{\infty}$-norm of controlled vs uncontrolled solutions relative to $\bar{h}$.}
\label{wave_norm_cost}
\end{figure}

We also carried out simulations aiming to retain the surface at a linear state, while solving the nonlinear system (recall that this exhibits free-surface instabilities in the uncontrolled case if $L>2\pi$).
The solution of the uncontrolled forward system with initial conditions \eqref{IC} and $h_A=0.5$, $m=3$, \(L = \frac{15}{2} \pi \) results in nonlinear instability, but our objective is to maintain the free surface at the linear initial condition \(\mathcal{H}(x,0)=h_0+s_0\) (i.e., \(\bar{h} = \mathcal{H}(x,0) \)) throughout the simulation period \([0, T]\), with final simulation time $T=10$.
The results of the computation are depicted in Figures~\ref{snapshots_end}, \ref{fig:wave2} and \ref{wave_norm_cost}.
Panel~\ref{EX4_hbar} shows the target state \(\bar{h}\), whereas panel~\ref{wave_H} presents $\mathcal{H}^*$ maintaining its initial sinusoidal shape. Panel~\ref{wave_h} depicts the free surface evolution, which naturally tends towards instability (non-flat state). In contrast, it is adjusted by the changing topography in panel~\ref{wave_s}, which develops inverse undulations to \(h^*\) to maintain $\mathcal{H}^*$. The oscillatory behaviour of the control $f^*$ is depicted in panel~\ref{wave_f}. 
The cost function decreases as indicated in panel~\ref{wave_cost}, reflecting progress towards the target state. Finally, Figure~\ref{wave_norm_cost} shows a brief increase in the infinity norm due to strong instability, before it returns to near zero. 

%%\clearpage

  \subsection{Rupture scenarios}
\label{rupture}

In this section, we set $\beta = 0$ in \eqref{reduced_cost}, focusing on controlling $h(x,t)$ rather than the total surface $\mathcal{H} = h + s$. This approach lets $\mathcal{H}$ rupture naturally, driven by the interaction between the thin film $h(x,t)$ and the underlying topography $s(x,t)$, providing a realistic representation of physical rupture. If we were to control $\mathcal{H}$ directly, it could prevent rupture by limiting the film’s ability to interact with the topography closely.

In the following examples, we aim to examine how thin-film flows evolve towards a rupture state. 
We set a non-zero rupture function $\phi(h)$ as in \eqref{IMEX}, with $\phi_{+} = \frac{1}{2 \varepsilon^4} h^2$, $\varepsilon = 0.1$, and $A = 0.03$. We integrate the forward system \eqref{Mixed} with $f=s=0$ and using the initial conditions \eqref{IC} with $h_A=0.5$, and within a spatial domain of length $L = 3\pi$. 
The simulations are carried out over $[0, T]$, with $T = 550$, until the dynamics settle into a steady configuration following rupture. Here, rupture refers to the film thickness falling below the regularisation threshold imposed through the regularised rupture function $\phi(h)$ as in~\eqref{regul0}; the parameter $\varepsilon$ is chosen to be small so that the dynamics towards the moment of film rupture are accurately captured, after which the evolution continues towards a steady configuration.
This numerically obtained quasi-steady profile is then used as the target state $\bar{h}_{\text{rupt}}$ for the corresponding optimal-control problem.

\begin{figure}[!p]
  \centering
  \hspace{-1.5cm}
  \begin{subfigure}[b]{0.4\textwidth}
    \centering
  \includegraphics[width=0.8\linewidth]{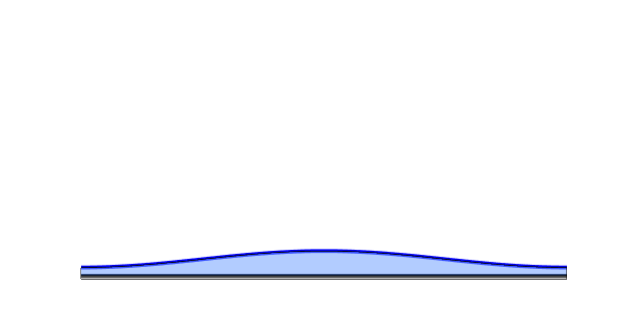}%{images2/filled_ruptureT70_HS_T0.pdf}
      \caption{$t=0$}
  \end{subfigure}
  \hspace{-1.5cm}  
  \begin{subfigure}[b]{0.4\textwidth}
    \centering
    \includegraphics[width=0.8\linewidth]{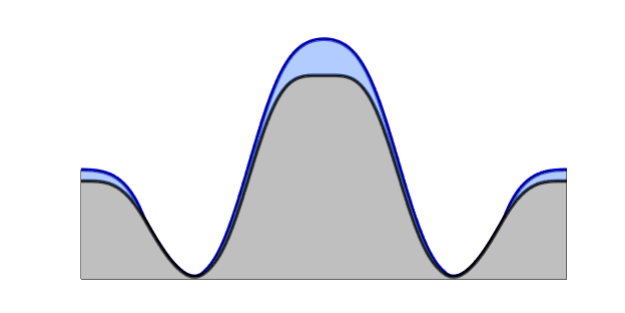}%{images2/filled_ruptureT70_HS_T18.pdf}
    \caption{$t=15$}
  \end{subfigure}
  \hspace{-1.5cm}  
  \begin{subfigure}[b]{0.4\textwidth}
    \centering
    \includegraphics[width=0.8\linewidth]{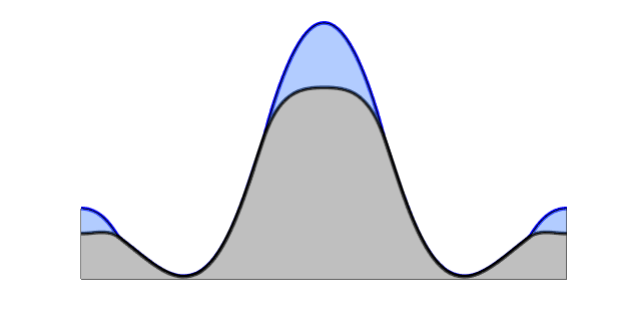}%{images2/filled_ruptureT70_HS_Tend.pdf}
    \caption{$t=30$}
  \end{subfigure}
  \caption{Example~6. Snapshots of the optimal solution, controlled towards a ruptured state. The solution is extended symmetrically in domain $[-L,L]$.}
  \label{snapshots_Hammondr111}
\end{figure}

\begin{figure}[!p]
  \centering
  \hspace{-1.0cm}
  \begin{subfigure}[b]{0.4\textwidth}
    \centering
\includegraphics[width=0.75\linewidth]{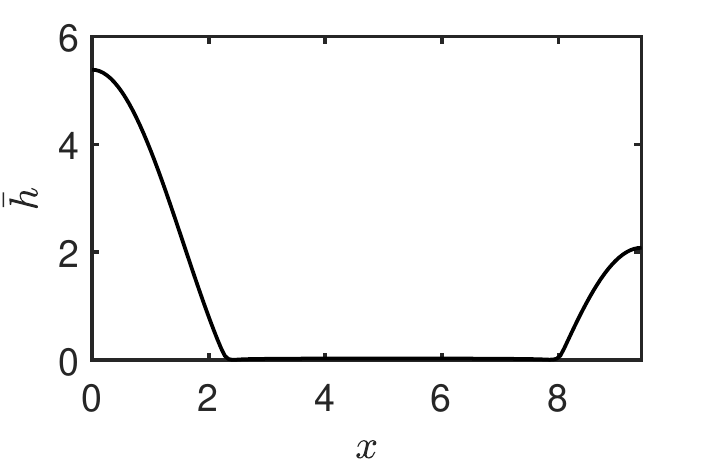}
    \caption{}
    \label{fig:sub1r_rupture}
  \end{subfigure}
  \hspace{-1.2cm}
  \begin{subfigure}[b]{0.4\textwidth}
    \centering  
    \includegraphics[width=0.7\linewidth]{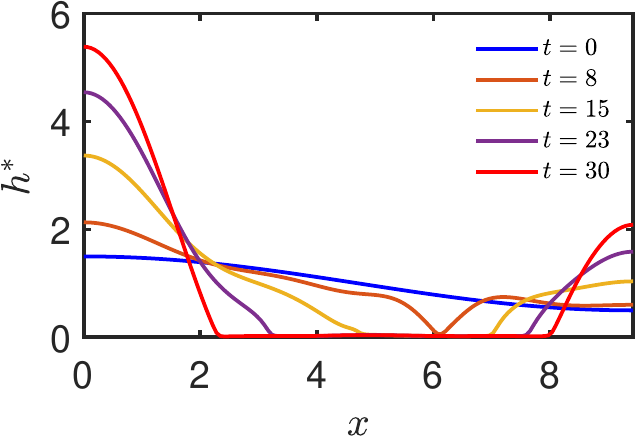}
    \caption{}
    \label{fig:sub4rr_rupture}
  \end{subfigure}
  \hspace{-1.2cm}
  \begin{subfigure}[b]{0.4\textwidth}
    \centering  
    \includegraphics[width=0.7\linewidth]{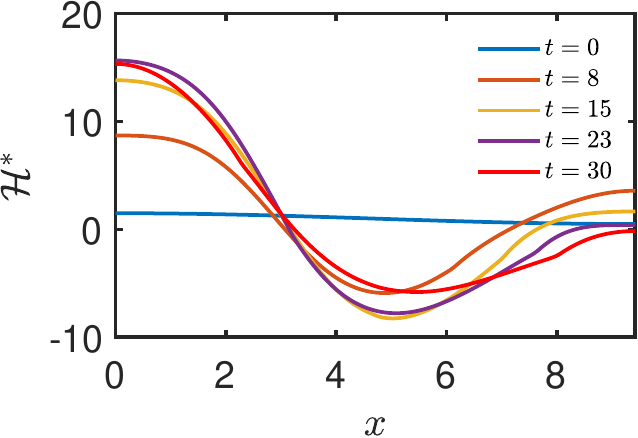}
    \caption{}
    \label{fig:sub4rr_ruptureHS}
  \end{subfigure}
 
   \hspace{-1.5cm}
  \begin{subfigure}[b]{0.4\textwidth}
    \centering
    \includegraphics[width=0.8\linewidth]{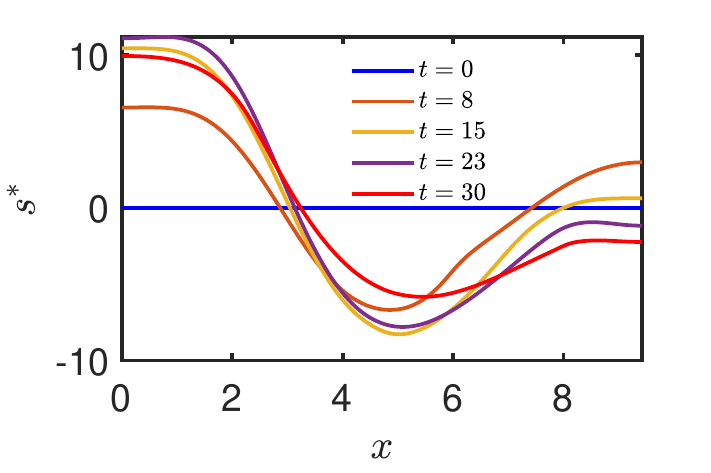}
    \caption{}
    \label{fig:sub4rrs}
  \end{subfigure}
  \hspace{-1.0cm}
    \begin{subfigure}[b]{0.4\textwidth}
    \centering
\includegraphics[width=0.8\linewidth]{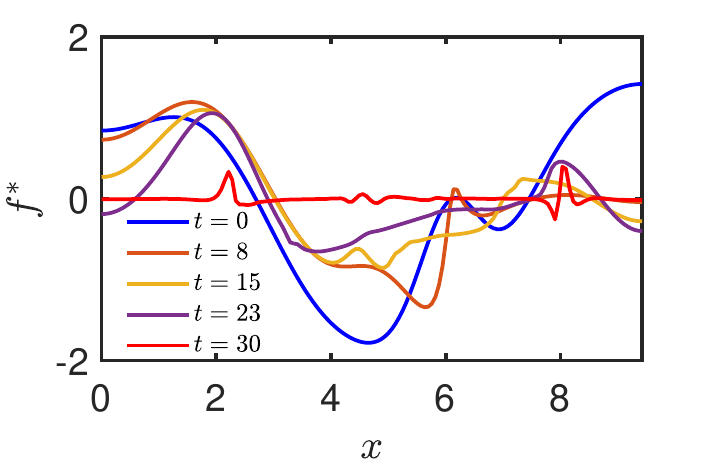}
    \caption{}
    \label{fig:hammond_fr_rupture}
  \end{subfigure}
 \caption{Example~6. (a) Target state $\bar{h}$. (b-e) Evolution of the optimal surface profile \( \mathcal{H}^* \), film height \( h^* \), topography \( s^* \), and control \( f^* \), respectively.}
  \label{rupture1}
\end{figure}

\subsubsection*{Example~6} 
\begin{figure}[!th]
\begin{subfigure}[b]{0.4\textwidth}
    \centering
    \includegraphics[width=0.8\linewidth]{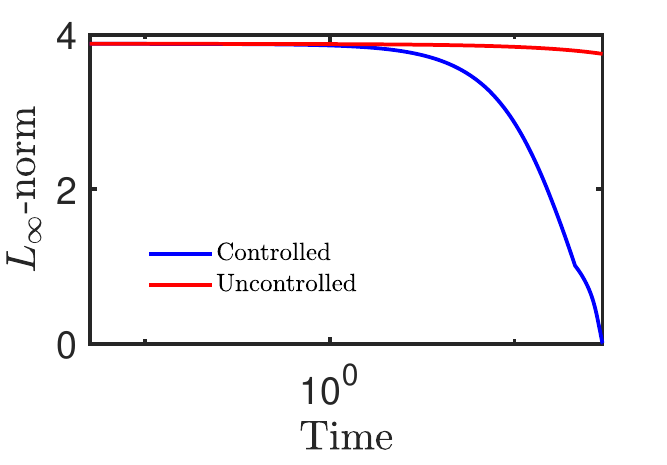}
    \caption{}
    \label{fig:sub5r}
  \end{subfigure}
  \centering
  \begin{subfigure}[b]{0.4\textwidth}
    \centering    
    \includegraphics[width=0.8\linewidth]{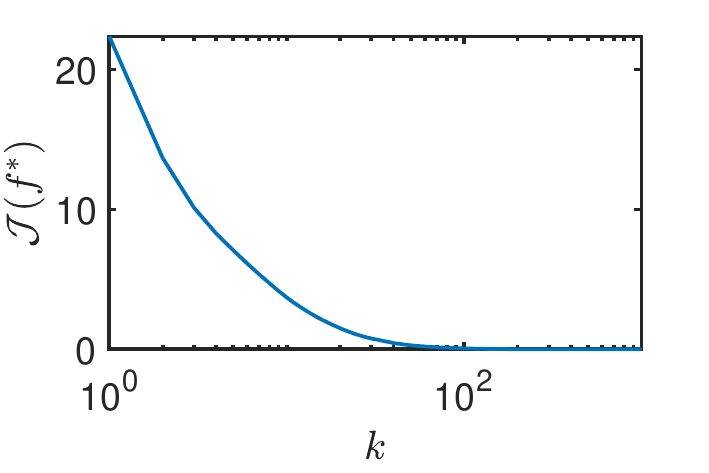}
    \caption{}
    \label{fig:sub6rr}
  \end{subfigure}
\caption{Example~6. (a) $L^{\infty}$-norm of the variation $h^*-\bar{h}_{\text{rupt}}$. (b) Cost function \(\mathcal{J}(f^*)\) against number of iterations \(k\).}
\label{cost_reupture}
\end{figure}

In the first example, we assess how quickly the optimal free surface $h^*$ reaches the target rupture state $\bar{h}=\bar{h}_{\text{rupt}}$, relative to the uncontrolled case. 
The optimal solution is shown in Figures~\ref{snapshots_Hammondr111} and \ref{rupture1}. 
While the uncontrolled solution necessitated $T=550$ to achieve a quasi-steady ruptured thin-film state $\bar{h}_{\text{rupt}}$, shown in panel~\ref{fig:sub1r_rupture}, our controlled solution $h^*$ achieves the desired state within just $T=30$, as seen in panel~\ref{fig:sub4rr_rupture}. 
Panel~\ref{fig:sub4rr_ruptureHS} depicts how the optimal thin-film flow over topography $\mathcal{H}^*$, including rupture effects, changes over time. 
The optimal topography $s^*$, and the control function $f^*$ are illustrated in panels \ref{fig:sub4rrs} and \ref{fig:hammond_fr_rupture}, respectively. 
The significant reduction in convergence time is highlighted in Figure~\ref{cost_reupture}. Both $L^{\infty}$ norm of $h^*-\bar{h}_{rupt}$ and cost function approach zero, though more slowly than in the non-rupture case.

\begin{figure}[!p]
  \centering
  \hspace{-1.5cm}
  \begin{subfigure}[b]{0.4\textwidth}
    \centering
  \includegraphics[width=0.8\linewidth]{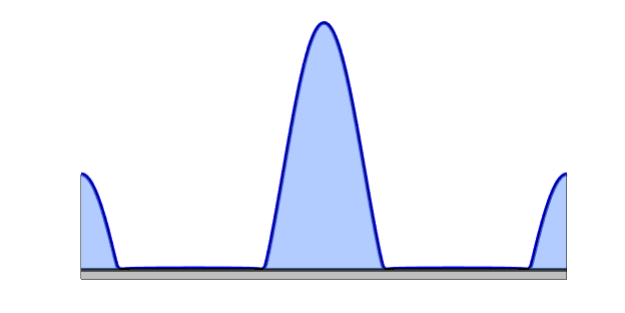}%
      \caption{$t=0$}
  \end{subfigure}
  \hspace{-1.5cm}  
  \begin{subfigure}[b]{0.4\textwidth}
    \centering
    \includegraphics[width=0.8\linewidth]{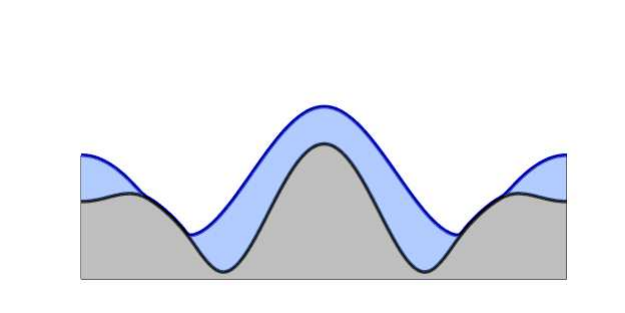}%
    \caption{$t=15$}
  \end{subfigure}
  \hspace{-1.5cm}  
  \begin{subfigure}[b]{0.4\textwidth}
    \centering
    \includegraphics[width=0.8\linewidth]{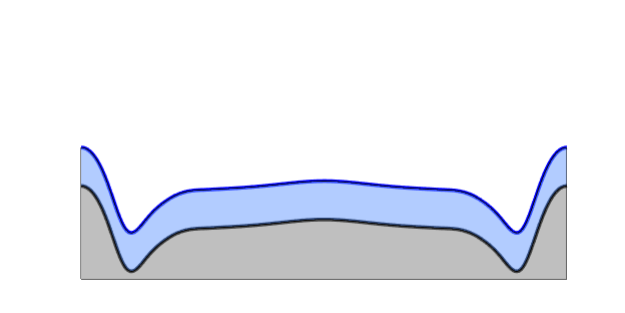}%
    \caption{$t=30$}
  \end{subfigure}
  \caption{Example~7. Snapshots of the optimal solution, controlled from a ruptured state to a uniform film thickness. The solution is extended symmetrically in domain $[-L,L]$.}
  \label{snapshots_Hammondr113}
\end{figure}

\begin{figure}[!p]
  \centering
  \hspace{-1.0cm}
  \begin{subfigure}[b]{0.4\textwidth}
    \centering  
    \includegraphics[width=0.75\linewidth]{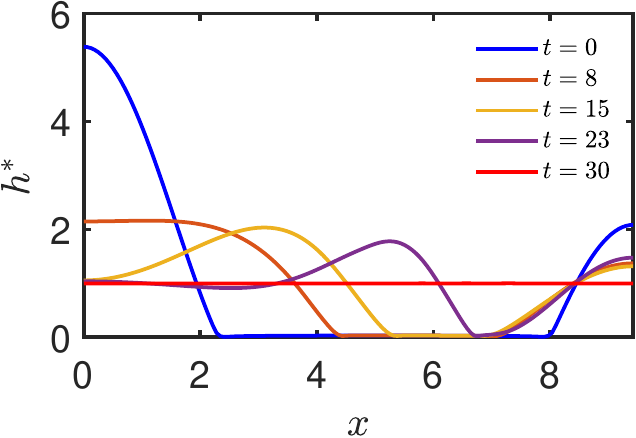}
    \caption{}
    \label{fig:sub_hr2}
  \end{subfigure}
  \hspace{-1.0cm}
   \begin{subfigure}[b]{0.4\textwidth}
   \centering 
\includegraphics[width=0.83\linewidth]{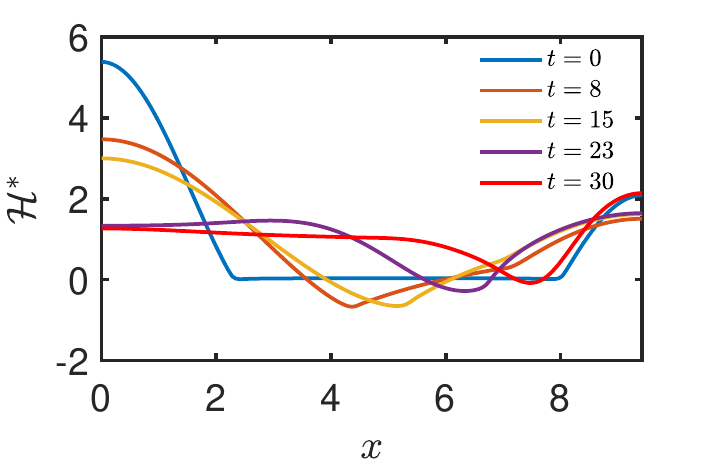}
    \caption{}
    \label{fig:sub33r}
  \end{subfigure}
  \hspace{-1.2cm}
  \begin{subfigure}[b]{0.4\textwidth}
    \centering
\includegraphics[width=0.8\linewidth]{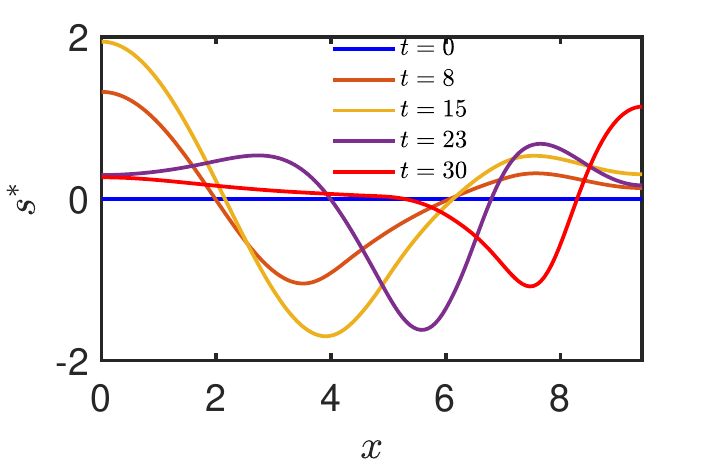}
    \caption{}
    \label{fig:sub4r2}
  \end{subfigure}

  \hspace{-1.0cm}
    \begin{subfigure}[b]{0.43\textwidth}
    \centering
\includegraphics[width=0.8\linewidth]{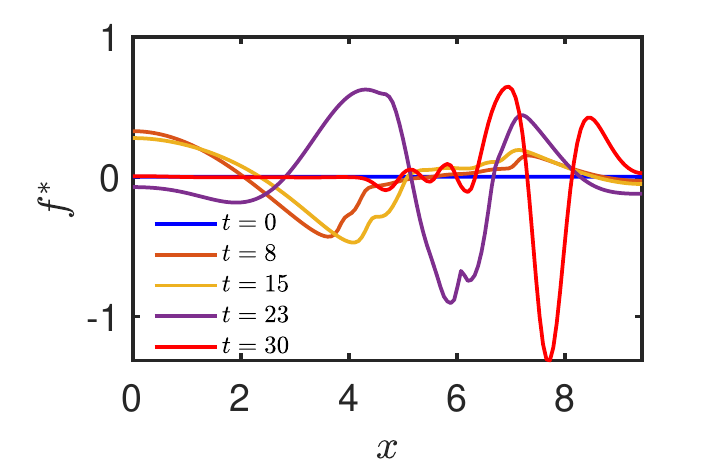}
    \caption{}
    \label{fig:hammond_fr2}
  \end{subfigure}
  \hspace{-1.2cm}
   \centering
  \begin{subfigure}[b]{0.4\textwidth}
    \centering    \includegraphics[width=0.8\linewidth]{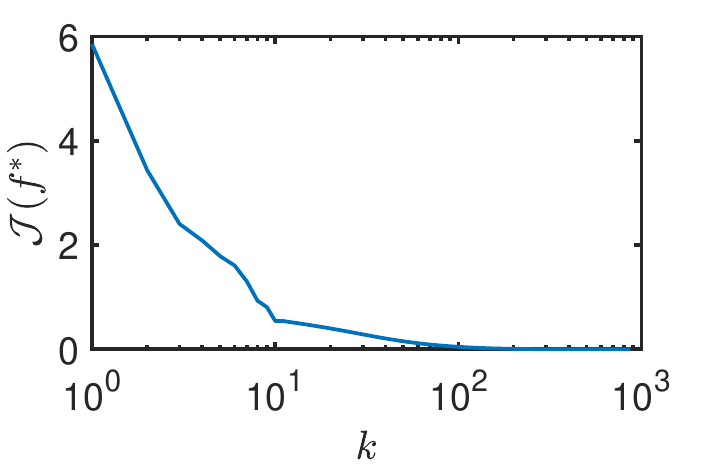}
    \caption{}
    \label{fig:sub6r2}
  \end{subfigure}
  \caption{Example~7. (a-d) Evolution of the optimal surface profile \( \mathcal{H}^* \), film height \( h^* \), topography \( s^* \), and control \( f^* \), respectively. (e) Cost function \(\mathcal{J}(f^*)\) against number of iterations \(k\).}
\label{flatting}
\end{figure}

\subsubsection*{Example~7} In the final example, we start with an initially ruptured state, $h_0(x)=\bar{h}_\text{rupt}$, and show that our control can halt the rupture by allowing the bubbles to merge into a single uniform film, as shown in Figure \ref{snapshots_Hammondr113}. This demonstrates how control can guide the system from rupture to a stable state. 
Figure~\ref{flatting} depicts $h^*$, $\mathcal{H}^*$, the flexible substrate $s^*$ and control $f^*$, as well as the cost function which is seen to reach its minimum faster compared to the previous example (cf. Figure~\ref{fig:sub6rr}).

 \section{Discussion}
 \label{Discussion}
 We have developed a mathematical framework for the open-loop optimal control of thin-film flows over flexible substrate, accounting for the rupture of films and the coalescence of bubbles. The focus of this study has been on finding the distributed force influencing the flexible substrate to minimise the cost function. We proposed a comprehensive model that couples thin-film and flexible-substrate deformations.  
We have also established a functional setting for the coupled PDE system, comprising a thin film equation and a flexible substrate equation. This demonstrates that this system adheres to a global energy-dissipation law. Furthermore, we have derived the optimality conditions by formulating the mathematical control problem as a PDE-constrained optimisation to obtain the time-dependent forward and backward PDE system and the optimality condition. We used a stable time discretisation approach with (IMEX) schemes to address the model’s nonlinearity, such as the disjoining-pressure term. This approach worked well in our numerical examples, demonstrating its effectiveness in solving the optimisation problem.

The effectiveness of our optimal control problem was demonstrated through a range of numerical examples, showing its ability to accelerate convergence to steady states, stabilise dewetting processes, achieve specific profiles, and even reverse rupture dynamics to restore uniformity. 
These results were achieved by allowing the adjacent wall to take a flexible shape, thereby affecting the surface profile. 
This non-standard approach provides a mechanism to control thin-film dynamics without altering the fluid properties (e.g., surface tension, stresses). 
It offers a useful alternative for addressing challenges in thin-film control and complements existing control approaches; previous studies have focussed on using thermal feedback \cite{garnier2003optical,Thompson_heat_2019}, fluid injection and extraction \cite{lou2003optimal,cimpeanu2021active,gomes2017stabilizing,thompson2016stabilising}, or electric fields \cite{tomlin2020instability,wray2022electrostatic}.

%While these approaches have made progress, they target the fluid properties directly. In this work, we took a different approach by focusing on how a flexible topography can be used to control thin-film behaviour.

The optimality system presented in this study was formulated as an open-loop control problem, namely the controller design is independent of the system output.
%the final result depends on the initial state.
This work can be extended to closed-loop control, also known as feedback control, which influences the control input directly~\cite{azmi2021optimal} and helps navigate the system along the optimal path set by the open-loop control problem~\cite{lunz2021minimizing}.
In thin film applications, feedback control not only stabilises contact line instabilities but also addresses complex problems, such as optimising the placement of actuators and sensors to suppress or enhance waves on a thin film~\cite{garnier2003optical,gomes2017stabilizing,thompson2016stabilising}. In addition, feedback control using fluid injection and extraction techniques are employed to mitigate wave phenomena in thin films and optimise actuator placement~\cite{holroyd2024linear,thompson2016stabilising}. 
The extension of our work to feedback control and developing numerical methods for solving a closed-loop control problem is left as future work.

 %\clearpage

\appendix

%\addtocontents{toc}{Appendix}

%\chapter{Derivation of flexible wall model}
%\input{appendices/Wall_Derive}

\section{Derivation of the thin film flow on a flexible substrate}
\label{app:derivation}
The derivation follows lubrication-theory for thin films on flexible substrates, as in Matar \emph{et al.}~\cite{matar2007dynamics}, adapted to the present configuration. We consider a Newtonian liquid film of typical thickness \(h_{0}\)
flowing beneath a flexible substrate over a horizontal length scale
\(L = h_{0}/\varepsilon\), where
\(
\varepsilon \coloneq \frac{h_{0}}{L} \ll 1
\)
is the small aspect ratio. 
The coordinate system is defined by \((x,y) \), with $x$ being the horizontal coordinate 
and $y$ being the vertical coordinate.
We substitute the scalings from Table~\ref{tab:nondim2} into the non-dimensional Navier-Stokes equations, the tangential and normal stress balance conditions, and the substrate deformation equation, in which an external forcing is included as part of the actuation pressure acting on the wall. 

Expanding the equations in powers of \(\epsilon\) and retaining leading-order terms yields:
\begin{subequations}
\begin{equation}\label{eq:p_x}
  p_{x} = u_{yy},
\end{equation}
\begin{equation}\label{eq:p_y}
  p_{y} = \frac{B_{o}}{C_a},
\end{equation}
\begin{equation}\label{eq:v_st}
  u_{x} + v_{y} = 0,
\end{equation}
\end{subequations}
where $p(x,y,t)$ and $\boldsymbol{u}=(u(x,y,t),\,v(x,y,t))$ are, respectively, the pressure and velocity of the fluid. All dimensionless quantities in this and subsequent equations ($C_a$, $B_o$, etc.) are defined in Table~\ref{tab:parameters}.
The leading-order tangential stress and the normal stress balance conditions at the free surface \(y=\mathcal{H}(x,t)=h(x,t)+s(x,t)\) are  
\begin{subequations}
\begin{align}
  u_y &= 0,   \label{eq:surface_bc} \\
  p &= -\Pi(h) - \frac{1}{C_a}\,\mathcal{H}_{xx},\label{eq:pressure}
\end{align}
\end{subequations}
 where the dimensionless disjoining
pressure is written in terms of the potential \(\phi\) as
\begin{equation} \label{Pi}
    \Pi(h) = -\phi^{\prime}(h) =Ah^{-3}.
\end{equation}  
The conditions on the wall~\(y=s(x,t)\) are
\begin{subequations}
\begin{align}
  u &= 0, \quad v = s_{t},  \label{eq:wall_bc} \\
  s_t - c^2 s_{xx} &= f(x,t) - \gamma p(x,s,t),
  \label{eq:substrate}
\end{align}
\end{subequations}
where the forcing term \(f(x,t)\) is defined as \(f = p_{\text{ext}} \). 

\begin{table}[t]
\centering
\caption{Variable scaling (asterisks denote dimensional quantities). \\ Here, $\bar{\mu}$ is the viscosity and $U_0$ is the characteristic velocity used in the scaling.}
\label{tab:nondim2}
\begin{tabular}{|>{\columncolor{white}}c|c|}
\hline
Dimensional variable & Nondimensional scaling \\[3pt]
\hline
$x^*$ & $x^* = (h_0/\varepsilon)\,x$ \\[3pt]
$y^*$ & $y^* = h_0\,y$ \\[3pt]
$u^*$ & $u^* = U_0\,u$ \\[3pt]
$v^*$ & $v^* = \varepsilon U_0\,v$ \\[3pt]
$t^*$ & $t^* = h_0/(\varepsilon U_0)\,t$ \\[3pt]
$p^*$ & $p^* = (\bar\mu U_0/(\varepsilon h_0))\,p$ \\[3pt]
\(\Pi^*\) & \(\Pi^* = (\bar\mu U_0/(\varepsilon h_0))\,\Pi\) \\[3pt]
$p_{ext}^*$ & $\displaystyle p_{ext}^*=\frac{p_{\mathrm{ext}}}{\bar\mu U_0/(\varepsilon h_0)}$ \\[3pt]
\hline
\end{tabular}
\end{table}

The kinematic condition
\begin{equation}
\mathcal{H}_t + u \mathcal{H}_x = v,
\end{equation}
can be written in the form 
\begin{equation} \label{flux_t}
    h_t + q_x = 0,
\end{equation}
where \(q\) is the flux defined by 
\begin{equation} \label{flux}
q(x,t) = \int_s^{s+h} u\,dy.
\end{equation}
To find $u$, we integrate \eqref{eq:p_y} and use the condition \eqref{eq:pressure}, which gives
\begin{equation}
p(x,y,t)
= -\Pi(h) - \frac{1}{C_a} \mathcal{H}_{xx}
  + \frac{B_o}{C_a}(y-\mathcal{H}).
\label{eq:p_field}
\end{equation}
Differentiating in $x$ yields
\begin{equation}
p_x
= -\Bigl[\Pi(h) + \frac{1}{C_a} \mathcal{H}_{xx}
             + \frac{B_o}{C_a} \mathcal{H}\Bigr]_x,
\label{eq:px_total}
\end{equation}
which shows that $p_x$ is independent of $y$ and hence we can simply integrate equation~\eqref{eq:p_x} and use the boundary condition~\eqref{eq:surface_bc}. Then we integrate the result once again and use the no-slip boundary condition~\eqref{eq:wall_bc} to get the final result
\begin{equation}
u(x,y,t) = p_x\,(y-s)\Bigl(\tfrac12 (y-s) - h\Bigr).
\end{equation}
Substituting the velocity into equation~\eqref{flux} gives the flux 
\begin{equation}
q(x,t) = \int_s^{s+h} u\,dy = -\frac{h^3}{3}\,p_x.
\end{equation}
Using equation~\eqref{eq:px_total} and equation~\eqref{flux_t}, we obtain
\begin{equation}
h_t + \frac{1}{3}\bigl(h^3 \Bigl[\Pi(h) + \frac{1}{C_a} \mathcal{H}_{xx}
             + \frac{B_o}{C_a} \mathcal{H}\Bigr]_x\bigr)_x = 0,
\label{eq:ht_G}
\end{equation}
which is the same as~\eqref{ht}, and~\eqref{mu1} is defined as the term in the square brackets.

To complete the coupling with the deformable substrate, we evaluate the pressure field \eqref{eq:p_field} at the wall position $y=s$, which gives 
\begin{equation}
p(x,s,t)
= -\Pi(h) - \frac{1}{C_a} \mathcal{H}_{xx} - \frac{B_o}{C_a} h.
\end{equation}
Substituting this expression into the wall equation \eqref{eq:substrate} yields the corresponding evolution equation for the substrate deflection. In this step we neglect the term $\gamma\Pi(h)$, which is asymptotically small compared with the hydrodynamic loading. The resulting relation coincides with the flexible-wall model presented in \eqref{flexible} of \eqref{Mixed}. 
\begin{table}[t!]
\centering
\caption{Nondimensional parameters. Here, $\rho$ is the density, $g$ is the gravity, $\sigma$ is the surface tension, $A^*$ is the dimensional Hamaker constant, $\gamma_w$ is the wall damping, and $T_w$ is the wall tension.}
\label{tab:parameters}
\begin{tabular}{c c}
\hline
Parameter & Definition \\[3pt]
\hline
Capillary  & $\displaystyle C_a=\frac{\bar\mu U_0}{\sigma\varepsilon^3}$ \\[6pt]
Bond & $\displaystyle B_o=\frac{\rho g h_0^2}{\sigma\varepsilon^2}$ \\[6pt]
Hamaker   & $\displaystyle A=\frac{A^* \varepsilon h_0^2}{\bar\mu U_0}$ \\[6pt]
Wall tension & $\displaystyle c^2=\frac{\sigma+T_w}{\gamma_w U_0 L}$ \\[6pt]
Film damping & $\displaystyle \gamma=\frac{\bar\mu }{\varepsilon^2\gamma_w^* h_0}$  \\[3pt]
\hline
\end{tabular}
\end{table}
\section{Free energy dissipation} \label{energyApp}
Here we derive the identity \eqref{energy}. We recall the energy
\begin{equation}
\mathcal{E}(h,s)=\int_\Omega\left[\phi(h)-\frac{B_o}{2 C_a} h^2+\frac{1}{2 C_a} h_x^2+\frac{c^2}{2 \gamma} s_x^2+ \frac{1}{C_a}(h_xs_x - B_o hs)\right] \mathrm{d}x,
\end{equation}
The molecular energy term $\phi(h)$ represents forces between molecules in the liquid film and the substrate, such as van der Waals forces. The film naturally moves toward thicknesses where $\phi(h)$ is smallest. In typical rupture examples, one sets $\phi'(h) = Ah^{-3}$ with $A > 0$, which represents attractive van der Waals forces, or a suitable regularisation thereof, e.g.~\eqref{regul0}. For such choices, $\phi(h)$ tends to a minimum as $h \to 0$, meaning the energy becomes lower when the film gets thinner. This makes thinning energetically favourable, so rupture can occur. The hydrostatic term $-\frac{B_{o}}{2C_a}h^2$ represents the effect of gravity. Since the film hangs beneath the substrate and this term is negative, the energy decreases when the film becomes thicker. Therefore, when $B_{o} > 0$, gravity encourages the film to form thicker bulges rather than staying uniform. The capillary term $\frac{1}{2C_a}h_x^2$ penalises sharp slopes in the film surface, favouring smoother and flatter shapes. Finally, there is the fluid-wall interaction energy term $\frac{1}{C_a}(h_x s_x - B_{o} h s)$: the cross-slope term $h_x s_x$ reduces the energy when the film and wall tilt in the same direction, encouraging them to deform together; the hydrostatic coupling $-B_{o} h s$ is reduced when thicker parts of the film sit above downward bends in the wall, representing the weight of the film pushing down on the substrate.

The goal is to derive the rate of change of this energy in time
\begin{equation}
\frac{\mathrm{d} \mathcal{E}}{\mathrm{d} t} = \frac{\mathrm{d}}{\mathrm{d} t} \int_\Omega \left[\phi(h) + \frac{1}{2 C_a} h_x^2 + \frac{c^2}{2 \gamma} s_x^2 + \frac{1}{C_a} h_x s_x - \frac{B_o}{2 C_a} h^2 - \frac{B_o}{C_a} h s \right] \mathrm{d}x.
\end{equation}
We apply the chain rule and integrate by parts on the terms involving spatial derivatives, then apply the boundary conditions \eqref{BC@} to remove the boundary terms, yielding
\begin{equation}
\frac{\mathrm{d} \mathcal{E}}{\mathrm{d} t} = \int_\Omega \left[ \left(\phi^{\prime}(h)-\frac{1}{C_a} h_{x x}-\frac{B_o}{C_a} h \right)   h_t - \frac{c^2}{\gamma} s_{xx} s_{t} + \frac{1}{C_a} (h_x s_x)_t - \frac{B_o}{C_a} (h s)_t  \right] \mathrm{d}x.
\end{equation}
We use equations \eqref{Mixed} to simplify the above expression, resulting in the dissipation law
\begingroup
\small
\begin{equation} \label{energy22}
\begin{aligned}
  \frac{\mathrm{d}}{\mathrm{d} t} \mathcal{E}(h,s)=  -\left(\frac{1}{3} \int_\Omega h^3 \mu_x^2 \,\mathrm{d}x + \frac{1}{\gamma} \int_\Omega s_t^2 \,\mathrm{d}x\right)+ \frac{1}{\gamma} \int_\Omega f s_t \,\mathrm{d}x.
\end{aligned}
\end{equation}
\endgroup
This equation represents the rate of change of the free energy. The first term on the right-had side %$\frac{1}{3} \int_\Omega h^3 \mu_x^2 \,\mathrm{d}x$ 
corresponds to the dissipation in the fluid, the second term %$ \frac{1}{\gamma} \int_\Omega s_t^2 \, \mathrm{d} x$ 
represents dissipation in the substrate, and the last term %$\frac{1}{\gamma} \int_\Omega f s_t \, \mathrm{d} x$ 
accounts for the work done by the external force on the substrate. Thus, without external forces (i.e., $f = 0$), the total free energy is non-increasing in time.

\section{Boundedness of \(( \phi'(h), w )\)}

We prove that for any \(h,w \in L^{2}(\Omega)\),
\begin{equation} 
    \label{ineq:phi_prime_bound}
    |( \phi'(h), w )|
    \le C_{\mathrm{L}} \bigl( \|h\|_{L^{2}(\Omega)} + C \bigr)\, \|w\|_{L^{2}(\Omega)},
\end{equation}
where \(C_{\mathrm{L}}\) is the Lipschitz constant of \(\phi'\) and \(C\) is a positive constant.

 Let \(h_\star\in\mathbb{R}\) be any fixed reference value. Recall that \(\phi^{\prime}\) is assumed to be  Lipschitz continuous. Hence, the following decomposition is valid:
\begin{equation}
    (\phi'(h), w)
    = (\phi'(h) - \phi'(h_\star),\, w)
      + (\phi'(h_\star), w).
\end{equation}

Using Lipschitz continuity of \(\phi'\),
\begin{align}
    | \phi'(h(x)) - \phi'(h_\star) |
    &\le C_{\mathrm{L}}\, |h(x) - h_\star|,
\end{align}
and applying Hölder's inequality gives
\begin{align}
    |(\phi'(h) - \phi'(h_\star), w)|
    &\le C_{\mathrm{L}} \|h - h_\star\|_{L^{2}(\Omega)} \, \|w\|_{L^{2}(\Omega)}.
\end{align}

Since \(h_\star\) is a constant,
\begin{align}
    |(\phi'(h_\star), w)|
    = |\phi'(h_\star)| \, \|w\|_{L^{1}(\Omega)}
    \le |\phi'(h_\star)|\, \|1\|_{L^{2}(\Omega)} \, \|w\|_{L^{2}(\Omega)}.
\end{align}

Combining both bounds,
\begin{align}
    |(\phi'(h), w)|
    &\le \bigl( C_{\mathrm{L}} \|h - h_\star\|_{L^{2}(\Omega)}
           + |\phi'(h_\star)|\,\|1\|_{L^{2}(\Omega)} \bigr)
           \|w\|_{L^{2}(\Omega)} \\
    &\le \bigl( C_{\mathrm{L}} (\|h\|_{L^{2}} + |h_\star|\, \|1\|_{L^{2}(\Omega)})
         + C^{\ast} \|1\|_{L^{2}(\Omega)} \bigr)
         \|w\|_{L^{2}(\Omega)}. 
\end{align}

Given that \(\|1\|_{L^{2}(\Omega)}= \sqrt{|\Omega|}\),  the inequality \eqref{ineq:phi_prime_bound} holds with
\begin{equation}
    C = \Big(C_{\mathrm{L}} |h_\star| + C^{\ast}\Big)\sqrt{|\Omega|},
\end{equation}
where \(h_\star\) may be chosen arbitrarily.
 \label{proof_in}
\section{IMEX-FEM discretisation} \label{FEM}

We discretise the spatial domain $\Omega \subset \mathbb{R}$ into $n-1$ elements
with $n$ nodes. The variational space $V=H^{1}(\Omega)$ is approximated by the
finite-dimensional subspace $V_{h}\subset V$ spanned by the standard
piecewise-linear hat functions $\{\hat{\phi}_i\}_{i=1}^{n}$.  
The solutions can therefore be characterised as linear combinations of $\hat{\phi}_{i}$ with real coefficients, such that $ \mathbf{F}_h(x,t)=\sum_{i=1}^{n} F_i(t)\,\hat{\phi}_i(x)$.
Here, $\mathbf{F}$ represents a vector that holds the unknowns at each node and is one of the state, adjoint or control variables, i.e. $\mathbf{F}=\{\mathbf{h},\boldsymbol{\mu},\mathbf{s},\mathbf{p},\mathbf{q},\mathbf{r},\mathbf{f}\}$, with corresponding coefficients $F_{i}=\{h_i,\mu_i,s_i,p_i,q_i,r_i,f_i\}$.

We define the real $n \times n$ matrices
\begin{align*}
&\mathbf{M}_{ij}=\!\int_{\Omega}\!\hat{\phi}_i\hat{\phi}_j\,dx,
\qquad
\mathbf{K}_{ij}=\!\int_{\Omega}\!\hat{\phi}_i'\hat{\phi}_j'\,dx,\\[2mm]
&\widehat{\mathbf{M}}_{ij}^{k}
   =\!\int_{\Omega}\!m(h_h^{k})\,\hat{\phi}_i\hat{\phi}_j\,dx, 
\qquad
\widehat{\mathbf{K}}_{ij}^{k}
   =\!\int_{\Omega}\!m(h_h^{k})\,\hat{\phi}_i'\hat{\phi}_j'\,dx,\\[2mm]
&\widehat{\mathbf{C}}_{ij}^{k}
   =\!\int_{\Omega}\!m(h_h^{k})\,\hat{\phi}_i'\hat{\phi}_j\,dx .
\end{align*}
The IMEX time discretisation uses uniform time steps
$t^{k}=k\Delta t$.
All linear terms are evaluated implicitly at $t^{k+1}$,
while nonlinear terms are evaluated explicitly at $t^{k}$.
In particular,~\(
\phi'_+(h_h^{k+1})\ \text{(implicit)},\ \text{and} \
\phi'_-(h_h^{k})\ \text{(explicit)},
\) and likewise for $\phi''_{\pm}$.

We write the semi-discrete optimality system in the following form:

\begin{subequations}\label{eq:IMEX-system}
\small
\begin{align}
\textbf{(State)}\qquad 
&\mathbf{M}\frac{\mathbf{h}^{k+1}-\mathbf{h}^{k}}{\Delta t}
   +\frac{1}{3}\widehat{\mathbf{K}}^{k}\bigl((\mathbf{h}^{k})^{3}\bigr)\boldsymbol{\mu}^{k+1}=0,
\\[1mm]
&\mathbf{M}\boldsymbol{\mu}^{k+1}
   -\frac{A}{\varepsilon^{4}}
     \widehat{\mathbf{M}}^{k+1}
        \bigl(\phi_{+}'(\mathbf{h}^{k+1})\bigr)\mathbf{h}^{k+1}
   -\frac{1}{C_{a}}\bigl(\mathbf{K}+B_{o}\mathbf{M}\bigr)\,
        \boldsymbol{\mathcal{H}}^{k+1}
\\[-1mm]
&\hspace{6cm}
=\frac{A}{\varepsilon^{4}}
     \widehat{\mathbf{M}}^{k}
       \bigl(\phi_{-}'(\mathbf{h}^{k})\bigr)\mathbf{h}^{k},
\\[1mm]
&\mathbf{M}\frac{\mathbf{s}^{k+1}-\mathbf{s}^{k}}{\Delta t}
   +c^{2}\mathbf{K}\mathbf{s}^{k+1}
   +\frac{\gamma}{C_{a}}\bigl(\mathbf{K}-B_{o}\mathbf{M}\bigr)\mathbf{h}^{k+1}
   =\mathbf{M}\mathbf{f}^{k+1},
\\[4mm]
\textbf{(Adjoint)}\qquad
&-\mathbf{M}\frac{\mathbf{p}^{k+1}-\mathbf{p}^{k}}{\Delta t}
   +\widehat{\mathbf{C}}^{k}\mathbf{p}^{k}
   +\left(
      \frac{B_o}{C_a}\mathbf{M}
      -\frac{1}{C_a}\mathbf{K}
      -\frac{A}{\varepsilon^{4}}
         \widehat{\mathbf{M}}^{k+1}
            \bigl(\phi_{+}''(\mathbf{h}^{k+1})\bigr)
     \right)\mathbf{q}^{k+1}
\\[-1mm]
&\qquad\qquad\qquad\qquad
+\frac{\gamma}{C_{a}}\bigl(B_{o}\mathbf{M}-\mathbf{K}\bigr)\mathbf{r}^{k+1}
 =\frac{A}{\varepsilon^{4}}
     \widehat{\mathbf{M}}^{k}
       \bigl(\phi_{-}''(\mathbf{h}^{k})\bigr)\mathbf{q}^{k}, 
\\[1mm]
&\frac{1}{3}\widehat{\mathbf{K}}^{k}\bigl((\mathbf{h}^{k})^{3}\bigr)\mathbf{p}^{k+1}
   +\mathbf{M}\mathbf{q}^{k+1}=0, \label{eq:disc_adj_q}
\\[1mm]
&-\mathbf{M}\frac{\mathbf{r}^{k+1}-\mathbf{r}^{k}}{\Delta t}
   +c^{2}\mathbf{K}\mathbf{r}^{k+1}
   -\frac{1}{C_{a}}\bigl(B_{o}\mathbf{M}-\mathbf{K}\bigr)\mathbf{q}^{k+1}
   =0,
\\[4mm]
\textbf{(Optimality)}\qquad
&\alpha\,\mathbf{M}\mathbf{f}^{k+1}
   -\mathbf{M}\mathbf{r}^{k+1}=0.
\end{align}
\end{subequations}

The forward problem is initialised with $\mathbf{h}^0$ and $\mathbf{s}^0$ from \eqref{IC}. The adjoint variables satisfy the terminal conditions obtained from the discrete optimality system, namely
\begin{equation}
\mathbf{p}^{N}=\mathbf{r}^{N}=\mathbf{h}^{N}+\beta\mathbf{s}^{N}-\bar{\mathbf{h}},
\end{equation}
with $\mathbf{q}^N$ given by \eqref{eq:disc_adj_q} at $k=N$.

\section*{Acknowledgments}
 SA acknowledges the financial support from Islamic University of Madinah and Ministry of Education of Saudi Arabia. The research by KvdZ was supported by the Engineering and Physical Sciences Research Council (EPSRC), UK under Grant and EP/W010011/1. We thank the anonymous reviewers for their careful reading and constructive comments, which helped to improve the clarity and quality of this work.

\bibliographystyle{siamplain}
\bibliography{references}

\begin{thebibliography}{10}

\bibitem{albi2017mean}
{\sc G.~Albi, Y.-P. Choi, M.~Fornasier, and D.~Kalise}, {\em Mean field control hierarchy}, Applied Mathematics \& Optimization, 76 (2017), pp.~93--135.

\bibitem{alexander2020stability}
{\sc J.~P. Alexander, T.~L. Kirk, and D.~T. Papageorgiou}, {\em Stability of falling liquid films on flexible substrates}, Journal of Fluid Mechanics, 900 (2020), p.~A40.

\bibitem{anderson2017electric}
{\sc T.~G. Anderson, R.~Cimpeanu, D.~T. Papageorgiou, and P.~G. Petropoulos}, {\em Electric field stabilization of viscous liquid layers coating the underside of a surface}, Physical Review Fluids, 2 (2017), p.~054001.

\bibitem{azmi2021optimal}
{\sc B.~Azmi, D.~Kalise, and K.~Kunisch}, {\em Optimal feedback law recovery by gradient-augmented sparse polynomial regression}, Journal of Machine Learning Research (JMLR), 22 (2021).

\bibitem{benney1966long}
{\sc D.~Benney}, {\em Long waves on liquid films}, Journal of mathematics and physics, 45 (1966), pp.~150--155.

\bibitem{bernis1990higher}
{\sc F.~Bernis and A.~Friedman}, {\em Higher order nonlinear degenerate parabolic equations}, Journal of differential equations, 83 (1990), pp.~179--206.

\bibitem{bertozzi1994lubrication}
{\sc A.~L. Bertozzi and M.~Pugh}, {\em The lubrication approximation for thin viscous films: the moving contact line with a'porous media'cut-off of van der waals interactions}, Nonlinearity, 7 (1994), p.~1535.

\bibitem{binard2022well}
{\sc J.~Binard, P.~Noble, and P.~Degond}, {\em Well-posedness and stability analysis of a landscape evolution model}, arXiv preprint arXiv:2211.09629,  (2022).

\bibitem{blyth2004effect}
{\sc M.~Blyth and C.~Pozrikidis}, {\em Effect of surfactant on the stability of film flow down an inclined plane}, Journal of Fluid Mechanics, 521 (2004), pp.~241--250.

\bibitem{blyth2023transition}
{\sc M.~G. Blyth, T.-S. Lin, and D.~Tseluiko}, {\em On the transition to dripping of an inverted liquid film}, Journal of Fluid Mechanics, 958 (2023), p.~A46.

\bibitem{camporeale2017asymptotic}
{\sc C.~Camporeale}, {\em An asymptotic approach to the crenulation instability}, Journal of Fluid Mechanics, 826 (2017), p.~636.

\bibitem{cimpeanu2021active}
{\sc R.~Cimpeanu, S.~N. Gomes, and D.~T. Papageorgiou}, {\em Active control of liquid film flows: beyond reduced-order models}, Nonlinear Dynamics, 104 (2021), pp.~267--287.

\bibitem{craster2009dynamics}
{\sc R.~V. Craster and O.~K. Matar}, {\em Dynamics and stability of thin liquid films}, Reviews of modern physics, 81 (2009), p.~1131.

\bibitem{de2015numerical}
{\sc J.~C. De~los Reyes}, {\em Numerical PDE-constrained optimization}, Springer, 2015.

\bibitem{decre2003gravity}
{\sc M.~M. Decr{\'e} and J.-C. Baret}, {\em Gravity-driven flows of viscous liquids over two-dimensional topographies}, Journal of Fluid Mechanics, 487 (2003), pp.~147--166.

\bibitem{duran2019instability}
{\sc M.~A. Dur{\'a}n-Olivencia, R.~S. Gvalani, S.~Kalliadasis, and G.~A. Pavliotis}, {\em Instability, rupture and fluctuations in thin liquid films: Theory and computations}, Journal of statistical physics, 174 (2019), pp.~579--604.

\bibitem{fowler2016controlling}
{\sc P.~D. Fowler, C.~Ruscher, J.~D. McGraw, J.~A. Forrest, and K.~Dalnoki-Veress}, {\em Controlling marangoni-induced instabilities in spin-cast polymer films: How to prepare uniform films}, The European Physical Journal E, 39 (2016), pp.~1--8.

\bibitem{frisk1972enhancement}
{\sc D.~P. Frisk and E.~J. Davis}, {\em The enhancement of heat transfer by waves in stratified gas-liquid flow}, International Journal of Heat and Mass Transfer, 15 (1972), pp.~1537--1552.

\bibitem{garnier2003optical}
{\sc N.~Garnier, R.~O. Grigoriev, and M.~F. Schatz}, {\em Optical manipulation of microscale fluid flow}, Physical Review Letters, 91 (2003), p.~054501.

\bibitem{gaskell_jimack_sellier_thompson_wilson_2004}
{\sc P.~Gaskell}, {\em Gravity-driven flow of continuous thin liquid films on non-porous substrates with topography}, Journal of Fluid Mechanics, 509 (2004), p.~253–280, \url{https://doi.org/10.1017/S0022112004009425}.

\bibitem{gaskell2004gravity}
{\sc P.~Gaskell, P.~Jimack, M.~Sellier, H.~Thompson, and M.~Wilson}, {\em Gravity-driven flow of continuous thin liquid films on non-porous substrates with topography}, Journal of Fluid Mechanics, 509 (2004), pp.~253--280.

\bibitem{gnann2018navier}
{\sc M.~V. Gnann and M.~Petrache}, {\em The navier-slip thin-film equation for 3d fluid films: existence and uniqueness}, Journal of Differential Equations, 265 (2018), pp.~5832--5958.

\bibitem{gomes2017stabilizing}
{\sc S.~N. Gomes, D.~T. Papageorgiou, and G.~A. Pavliotis}, {\em Stabilizing non-trivial solutions of the generalized kuramoto--sivashinsky equation using feedback and optimal control: Lighthill--thwaites prize}, IMA Journal of Applied Mathematics, 82 (2017), pp.~158--194.

\bibitem{GomZeeBOOK-CH2017}
{\sc H.~Gomez and K.~G. van~der Zee}, {\em Computational phase-field modeling}, in Encyclopedia of Computational Mechanics, Second Edition, E.~Stein, R.~de~Borst, and T.~J.~R. Hughes, eds., Wiley, 2017.
\newblock Part~1~Fluids.

\bibitem{halpern1992fluid}
{\sc D.~Halpern and J.~Grotberg}, {\em Fluid-elastic instabilities of liquid-lined flexible tubes}, Journal of Fluid Mechanics, 244 (1992), pp.~615--632.

\bibitem{hammond1983nonlinear}
{\sc P.~S. Hammond}, {\em Nonlinear adjustment of a thin annular film of viscous fluid surrounding a thread of another within a circular cylindrical pipe}, Journal of fluid Mechanics, 137 (1983), pp.~363--384.

\bibitem{herzog2010algorithms}
{\sc R.~Herzog and K.~Kunisch}, {\em Algorithms for pde-constrained optimization}, GAMM-Mitteilungen, 33 (2010), pp.~163--176.

\bibitem{holroyd2024linear}
{\sc O.~A. Holroyd, R.~Cimpeanu, and S.~N. Gomes}, {\em Linear quadratic regulation control for falling liquid films}, SIAM Journal on Applied Mathematics, 84 (2024), pp.~940--960.

\bibitem{hu2016effects}
{\sc B.~Hu}, {\em Effects of surface tension and viscoelastic behavior on the thin film coating flow of microbicide gels},  (2016).

\bibitem{jensen2004thin}
{\sc O.~E. Jensen, G.~Chini, and J.~King}, {\em Thin-film flows near isolated humps and interior corners}, Journal of engineering mathematics, 50 (2004), pp.~289--309.

\bibitem{kabaliuk2013blood}
{\sc N.~Kabaliuk, M.~Jermy, K.~Morison, T.~Stotesbury, M.~Taylor, and E.~Williams}, {\em Blood drop size in passive dripping from weapons}, Forensic science international, 228 (2013), pp.~75--82.

\bibitem{kalogirou2020nonlinear}
{\sc A.~Kalogirou and M.~G. Blyth}, {\em Nonlinear dynamics of two-layer channel flow with soluble surfactant below or above the critical micelle concentration}, Journal of Fluid Mechanics, 900 (2020).

\bibitem{kistler1997coating}
{\sc S.~F. Kistler and P.~M. Schweizer}, {\em Coating science and technology: An overview}, Liquid Film Coating,  (1997), pp.~3--15.

\bibitem{klein2016optimal}
{\sc M.~Klein and A.~Prohl}, {\em Optimal control for the thin film equation: Convergence of a multi-parameter approach to track state constraints avoiding degeneracies}, Computational Methods in Applied Mathematics, 16 (2016), pp.~685--702.

\bibitem{lee2011dynamics}
{\sc Y.~Lee, H.~Thompson, and P.~Gaskell}, {\em Dynamics of thin film flow on flexible substrate}, Chemical Engineering and Processing: Process Intensification, 50 (2011), pp.~525--530.

\bibitem{lou2003optimal}
{\sc Y.~Lou and P.~D. Christofides}, {\em Optimal actuator/sensor placement for nonlinear control of the kuramoto-sivashinsky equation}, IEEE Transactions on Control Systems Technology, 11 (2003), pp.~737--745.

\bibitem{lunz2021minimizing}
{\sc D.~Lunz}, {\em Minimizing deformation of a thin fluid film driven by fluxes of momentum and heat}, Physical Review E, 103 (2021), p.~033105.

\bibitem{manzoni2021optimal}
{\sc A.~Manzoni, A.~Quarteroni, and S.~Salsa}, {\em Optimal control of partial differential equations}, Springer, 2021.

\bibitem{matar2007falling}
{\sc O.~Matar, R.~Craster, and S.~Kumar}, {\em Falling films on flexible inclines}, Physical Review E, 76 (2007), p.~056301.

\bibitem{matar2004rupture}
{\sc O.~K. Matar and S.~Kumar}, {\em Rupture of a surfactant-covered thin liquid film on a flexible wall}, SIAM Journal on Applied Mathematics, 64 (2004), pp.~2144--2166.

\bibitem{matar2007dynamics}
{\sc O.~K. Matar and S.~Kumar}, {\em Dynamics and stability of flow down a flexible incline}, Journal of Engineering Mathematics, 57 (2007), pp.~145--158.

\bibitem{miles2020thermomechanically}
{\sc C.~Miles, K.~G. v.~d. Zee, M.~E. Hubbard, and R.~MacKenzie}, {\em Thermomechanically-consistent phase-field modeling of thin film flows}, Numerical Methods for Flows: FEF 2017 Selected Contributions,  (2020), pp.~121--129.

\bibitem{miyara1999numerical}
{\sc A.~Miyara}, {\em Numerical analysis on flow dynamics and heat transfer of falling liquid films with interfacial waves}, Heat and Mass Transfer, 35 (1999), pp.~298--306.

\bibitem{oron1997long}
{\sc A.~Oron, S.~H. Davis, and S.~G. Bankoff}, {\em Long-scale evolution of thin liquid films}, Reviews of modern physics, 69 (1997), p.~931.

\bibitem{reynolds1886iv}
{\sc O.~Reynolds}, {\em Iv. on the theory of lubrication and its application to mr. beauchamp tower’s experiments, including an experimental determination of the viscosity of olive oil}, Philosophical transactions of the Royal Society of London,  (1886), pp.~157--234.

\bibitem{rohlfs2017hydrodynamic}
{\sc W.~Rohlfs, P.~Pischke, and B.~Scheid}, {\em Hydrodynamic waves in films flowing under an inclined plane}, Physical Review Fluids, 2 (2017), p.~044003.

\bibitem{sellier2008substrate}
{\sc M.~Sellier}, {\em Substrate design or reconstruction from free surface data for thin film flows}, Physics of Fluids, 20 (2008).

\bibitem{sellier2016inverse}
{\sc M.~Sellier}, {\em Inverse problems in free surface flows: a review}, Acta Mechanica, 227 (2016), pp.~913--935.

\bibitem{sellier2010beating}
{\sc M.~Sellier and S.~Panda}, {\em Beating capillarity in thin film flows}, International journal for numerical methods in fluids, 63 (2010), pp.~431--448.

\bibitem{serifi2004transient}
{\sc K.~Serifi, N.~A. Malamataris, and V.~Bontozoglou}, {\em Transient flow and heat transfer phenomena in inclined wavy films}, International journal of thermal sciences, 43 (2004), pp.~761--767.

\bibitem{Thompson_heat_2019}
{\sc A.~B. Thompson, S.~N. Gomes, F.~Denner, M.~C. Dallaston, and S.~Kalliadasis}, {\em Robust low-dimensional modelling of falling liquid films subject to variable wall heating}, Journal of Fluid Mechanics, 877 (2019), p.~844–881, \url{https://doi.org/10.1017/jfm.2019.580}.

\bibitem{thompson2016stabilising}
{\sc A.~B. Thompson, S.~N. Gomes, G.~A. Pavliotis, and D.~T. Papageorgiou}, {\em Stabilising falling liquid film flows using feedback control}, Physics of Fluids, 28 (2016), p.~012107.

\bibitem{tomlin2020instability}
{\sc R.~J. Tomlin, R.~Cimpeanu, and D.~T. Papageorgiou}, {\em Instability and dripping of electrified liquid films flowing down inverted substrates}, Physical Review Fluids, 5 (2020), p.~013703.

\bibitem{tomlin2019optimal}
{\sc R.~J. Tomlin, S.~N. Gomes, G.~A. Pavliotis, and D.~T. Papageorgiou}, {\em Optimal control of thin liquid films and transverse mode effects}, SIAM Journal on Applied Dynamical Systems, 18 (2019), pp.~117--149.

\bibitem{Troeltzsch2010}
{\sc F.~Tr{\"o}ltzsch}, {\em Optimal Control of Partial Differential Equations: Theory, Methods and Applications}, vol.~112 of Graduate Studies in Mathematics, American Mathematical Society, Providence, RI, 2010.

\bibitem{tseluiko2013stability}
{\sc D.~Tseluiko, M.~G. Blyth, and D.~Papageorgiou}, {\em Stability of film flow over inclined topography based on a long-wave nonlinear model}, Journal of Fluid Mechanics, 729 (2013), pp.~638--671.

\bibitem{tseluiko2006wave}
{\sc D.~Tseluiko and D.~T. Papageorgiou}, {\em Wave evolution on electrified falling films}, Journal of Fluid Mechanics, 556 (2006), pp.~361--386.

\bibitem{wray2022electrostatic}
{\sc A.~W. Wray, R.~Cimpeanu, and S.~N. Gomes}, {\em Electrostatic control of the navier-stokes equations for thin films}, Physical Review Fluids, 7 (2022), p.~L122001.

\bibitem{wu2018posteriori}
{\sc X.~Wu, K.~G. van~der Zee, G.~Simsek, and E.~Van~Brummelen}, {\em A posteriori error estimation and adaptivity for nonlinear parabolic equations using imex-galerkin discretization of primal and dual equations}, SIAM Journal on Scientific Computing, 40 (2018), pp.~A3371--A3399.

\end{thebibliography}
\end{document}